\documentclass[smallextended,natbib,runningheads]{svjour3}

\usepackage{amsmath}
\usepackage{amssymb}
\usepackage{theorem}
\usepackage{xspace}
\usepackage[final]{graphicx}
\usepackage{subfigure}
\usepackage{xspace}
\usepackage{hhline}
\usepackage[a4paper]{geometry}
\theoremstyle{plain}

\usepackage{amsmath}
\usepackage{amssymb}
\usepackage{theorem}
\usepackage{xspace}
\usepackage{xspace}
\usepackage{hhline}
\usepackage{bbm}
\usepackage[a4paper]{geometry}

\def\pth#1{\left(#1\right)}
\def\acc#1{\left\{#1\right\}}
\def\cro#1{\left[#1\right]}

\def\eE{I\!\!E}
\def\eP{I\!\!P}
\def\e1{1\!\!1}

\def\eY{\textbf{Y}}
\def\eX{\textbf{X}}
\def\ef{\textrm{\mathversion{bold}$\mathbf{\phi}$\mathversion{normal}}}  
\def\ee1{\textrm{\mathversion{bold}$\mathbf{\varepsilon}$\mathversion{normal}}}  
\def\eth{\textrm{\mathversion{bold}$\mathbf{\theta}$\mathversion{normal}}}  
\def\eeps{\textrm{\mathversion{bold}$\mathbf{\varepsilon}$\mathversion{normal}}}

\def\eX{\mathbf{X}}
\def\ew{\mathbf{w}}
\def\eeX{\mathbb{X}}

\def\eO{\textrm{\mathversion{bold}$\mathbf{\Omega}$\mathversion{normal}}} 

\usepackage{theorem}
\usepackage{xspace}
\usepackage{xspace}
\usepackage{hhline}
\usepackage{bbm}
\usepackage[a4paper]{geometry}

\newcommand{\N}{\mathbb{N}}
\newcommand{\R}{\mathbb{R}}

\def\argmin{\mathop{\mathrm{arg\,min}}}

\begin{document}
\title {{\bf Model selection by LASSO methods in a change-point model}}

\author{GABRIELA CIUPERCA}
\institute{GABRIELA  CIUPERCA \at
Universit\'e de Lyon,
Universit\'e Lyon 1, 
CNRS, UMR 5208, Institut Camille Jordan, 
Bat.  Braconnier,
43, blvd du 11 novembre 1918,
F - 69622 Villeurbanne Cedex, France,\\\email{Gabriela.Ciuperca@univ-lyon1.fr}\\
{\it tel: }33(0)4.72.43.16.90, {\it fax: }33(0)4.72.43.16.87}
\date{Received: date / Accepted: date}
\maketitle


\begin{abstract}
The paper considers a linear  regression model with multiple change-points occurring at unknown times. 
The LASSO technique is very interesting since it allows the parametric estimation, including the change-points, and automatic variable selection simultaneously. The asymptotic properties of the LASSO-type (which has as particular case the LASSO estimator) and of the adaptive LASSO estimators are studied. For this last estimator the oracle properties are proved. In both cases, a model selection criterion is proposed. 
Numerical examples are provided showing the performances of the adaptive LASSO estimator compared to the LS estimator.
\end{abstract}


\section{Introduction}
A change-point model is a model which changes at the unknown observations. Change-point detection procedures fall into two categories: retrospective or a posteriori change detection and on-line, sequential or a priori change detection.  This paper focuses on an \textit{a posteriori} change-point problem, which arises when the data are completely known at the end of the experiment to process.
More precisely, we study the high-sized a posteriori change-point model:  study a phenomenon (dependent variable), function of one very large regressors variables number,  with unknown change-points number. \\
A significant advancement in variable selection in a model without change-point was realized by \cite{Tibshirani:96}, proposing the LASSO method. Then, the estimation and model selection are simultaneously treated as a single minimization problem. If the model have change-points, the LASSO method would allow at the same time  to estimate the parameters on every segment and eliminate the irrelevant predictive covariates without crossing every time by a hypothesis test. We remark that the least squares (LS) method gives nonzero estimates to all coefficients. \\
In this paper we propose a method for  a change-points linear model with the aim of estimating and choosing the covariates (regressors) simultaneously. The obtained results will be very useful for high-sized models which is often used in various fields especially in medicine, meteorology or financial econometrics. \\
Concerning the LASSO method in a model without change-points, a generalization of LASSO ($L_1$ penalty) and ridge ($L_2$ penalty) estimators was given in \cite{Knight:Fu:00} by minimizing the residual sum of squares plus a penalty proportion to the models parameters. The obtained estimator is called LASSO-type estimator. The ridge method has a good prediction performance through a bias-variance trade-off, while the LASSO method encourages both shrinkage and automatic variable selection simultaneously and it has very good computational properties. Nevertheless the LASSO estimator does not satisfy the oracle properties. To remedy this inconvenience, an adaptive LASSO  estimator   was proposed by \cite{Zou:06}. Recall that the oracle properties are: the zero components of the true parameters are estimated (shrunk) as 0 with probability tending to 1 (also called sparsity property) and the nonzero components have an optimal estimation rate (and is asymptotically normal). \\
Let us give some references of recent papers on the LASSO method, with the remark that to the author's knowledge, the LASSO problem has not yet been addressed in a change-point model. 
 For the models without change-points, in the case of linear regression, \cite{Fan:Li:01} show that the LASSO method produces biased estimator for the large parameter regression and \cite{Zou:06} proves that the oracle properties do not hold for the LASSO.  \cite{Potscher:Schneider:09} study the distribution of the adaptive LASSO estimator in finite samples and in the large-sample limit while \cite{Xu:Ying:10}  consider the LASSO-type penalty in a median regression. In the paper of \cite{Foster:Verbyla:Pitchford:09} a LASSO random effects model is considered and in \cite{Bickel:Ritov:Tsybakov:09}, the equivalence results and sparsity oracle inequalities for the LASSO and Dantzing estimators in a nonparametric regression model are given. \cite{Wei:Huang:Li:11} consider the problem of variable selection and estimation for a linear model with time-varying effects for covariates using the group LASSO and adaptive group LASSO methods. They proved that the obtained estimators are consistent and the adaptive group LASSO estimator has the oracle selection property. \\
 Seen the very interesting properties of these estimators (LASSO-type and adaptive LASSO), we study its behavior in a model with change-points. In order to estimate the change-point number, we also propose a model selection  criterion based on LASSO-type or adaptive LASSO method. In a multiple change-point model, the break estimation could affect the estimator properties: variable selection on each segment, oracle properties, ... This is the main interest of this paper. Besides, since the penalty contains the model parameters, the change-point results of the literature (see e.g. \cite{Bai:98}, \cite{Bai:Perron:98}, \cite{Koul:Qian:02} or \cite{Ciuperca:09}) do not apply. \\
This paper considers the estimation and the selection of the significant variables in a linear regression with multiple change-points occurring at unknown time. The study of this method was motivated by wishes to find the properties of the estimator, particularly interesting in a change-point model with high-sized regressors, which allows the automatic elimination of the non significant variables on every phase, without using hypothesis test.\\
Two estimation methods are proposed and studied: LASSO-type (with particular cases: ridge and LASSO method) and adaptive LASSO. The first method has the advantage that it can consider a lower observations size with respect to model parameter number. But, under certain conditions on design (matrix of observed regressors), the LASSO-type regression parameters estimators are not consistent. Then, we can use the adaptive LASSO method which correctly selects variables of nonzero coefficients with probability converging to one. On the other hand, in order to calculate an adaptive penalization, this method can be used  only if the observations size on every segment is bigger than the parameters number in corresponding interval. \\
The difficulty to study a change-point model, when the number of change-points is fixed, results first from the dependence of the model of two parameter types: the regression parameters and the change-points. Moreover, for multiple breaks, each middle regime has completely unknown boundaries. Then, for the two   estimators we study the asymptotic behavior: convergence rate and asymptotic distribution. For the adaptive LASSO regression parameter estimator, we also prove that the change-point presence in the model does not influence the oracle property: between two consecutive breaks, nonzero parameters estimator is asymptotically normal and zero parameters are shrunk directly to 0 with a probability converging to 1. \\
If the number of breaks is unknown, the problem of its estimate arises. We propose, a general criterion to estimate the number of change-points. \\
Finally, numerical simulations are realized to illustrate the theoretical results and to show the advantages of the proposed methods in terms of detection of irrelevant variables in a change-points model and also of break number estimation.\\  
In the present paper, as original contribution we provide statistical asymptotic properties of the LASSO-type and adaptive LASSO estimator in a change-point model.  The structure of this paper is as follows. The model and assumptions are introduced in Section 2. In Section 3, a LASSO-type estimator in a change-point model is proposed and its asymptotic behavior is studied. Convergence rate, asymptotic distribution of the regression parameters and of the change-point estimators are obtained. A model selection criterion is also studied.  Next, adaptive LASSO estimator and its oracle properties are given in Section 4. Section 5 reports some simulation results which illustrate the theoretical results. Finally, Appendix contains the proofs of the results.

\section{Model} 

We consider the  model: $Y_i=f_\eth(\eX_i)+\varepsilon_i$, for the step-function with $K$ ($K \geq 0$) change-points:
\[
f_\eth(\eX_i)=h_{\ef_1}(\eX_i) \e1_{i<l_1}+h_{\ef_2}(\eX_i) \e1_{l_1 \leq i<l_2}+ \cdots +h_{\ef_{K+1}}(\eX_i) \e1_{i>l_K}, \qquad i=1, \cdots, n
\]
more precisely $h_\ef(\eX)=\eX' \ef$, $\ef \in  \Gamma \subseteq \R^p $, $\Gamma$ compact.\\

Up to our knowledge, no result exists in the literature on the LASSO  estimation in a model with change-points. In all previous works, where the estimation by penalization is considered in a  multi-phase model, the penalization does not contain  model parameters (see \cite{Ciuperca:11b}) and thus we cannot make model selection and  the estimation simultaneously. Classically, in the rich literature on the change-point estimation, many articles contribute to determine the break number, to estimate the change-point location and the regression parameters (see e.g. \cite{Bai:Perron:98}, \cite{Kim:Kim:08} or \cite{Ciuperca:09}). Once the asymptotic distribution of estimator is proved, the hypothesis test may be performed to eliminate the irrelevant predictive covariates. This approach requires a huge amount of computation if $p$ or $K$  are large. A solution to  carry out change-point analysis, perform variable selection and estimate regression parameters simultaneously was proposed by \cite{Wu:08}. Since it is necessary to calculate all sub-models, of order $2^{p+4}$ when $K=1$, the criteria considered by \cite{Wu:08} still require many calculations. On the other hand, if $p$ is large compared to $n$, the Wu's criteria need to be modified.    \\
We can also recall the paper \cite{Harchaoui:Levy:10} where the estimation of the location of change-points in one-dimensional piecewise constant in an white noise is reformulated as a variable selection problem with a $L_1$ penalty. Note that their model is constant between two consecutive change-points and reformulated model is a classical linear regression without breaks. \\

 For the sample $i$, $Y_i$ denotes the response variable, $\eX_i$ is a p-vector of regressors and the $\varepsilon_i$ is the error. The errors $(\varepsilon_i)_{1 \leq i \leq n}$ are independent identically distributed(i.i.d.) random variables. The model parameters are $\eth\equiv(\ef_1, \cdots,\ef_{K+1},l_1,\cdots, l_K) \in \Gamma^{K+1} \times \mathbb{N}^K$ and their true values (unknown) are  $\eth^0 \equiv (\ef^0_1, \cdots, \ef^0_{K+1},$ $ l^0_1, \cdots, l^0_K)$.  
We set: $\eth\equiv(\eth_1,\eth_2)$, with $\eth_1\equiv(\ef_1, \cdots, \ef_{K+1})$ the regression parameters, $\eth_2\equiv(l_1, \cdots, l_K)$ the change-points.\\
Denote by $\phi_{r,k}$ the $k$th component of $\ef_r$ and $\phi^0_{r,k}$ the $k$th component of $\ef^0_r$, for $r \in \{ 1, \cdots, K+1 \}$, $k=1, \cdots, p$.\\
Let us consider the deterministic design matrix $\mathbb{X}\equiv(X_{ij})_{\overset{1\leq i \leq n}{1 \leq j \leq p}}$,  ${\eX^j}'$ its column $j$ and $\eX_i$ the $i$th line.\\
We now state the assumptions under which the asymptotic properties of estimators will be derived.\\

First, we impose the condition that the change-points are sufficiently far apart:\\
\textbf{(H1)} there exists two positive  constants $u, c(>0)$ such that $l_{r+1}-l_r \geq c_0 [n^u]$, for every $r=1,\cdots , K$,   with $l_0=1$ and $l_{K+1}=n$.\\
Without loss of generality, in the following we consider $3/4 \leq u \leq 1$ and $c_0=1$.\\

For the design  $\eeX$, we suppose that:\\
\textbf{(H2)} $n^{-1} \max_{1 \leq i \leq n} \eX_i' \eX_i {\underset{n \rightarrow \infty}{\longrightarrow}}  0$ and for any $r=1, \cdots, K+1$, the matrix \\
$\textbf{C}_{n,r}\equiv(l_r-l_{r-1} )^{-1} \sum^{l_r}_{i=l_{r-1}+1} \eX_i \eX_i' {\underset{n \rightarrow \infty}{\longrightarrow}} \textbf{C}_r$, with $\textbf{C}_r$ a non-negative definite matrix.\\

 For the errors $\varepsilon_i$ we suppose that:\\
 \textbf{(H3)} $\varepsilon$ is a random variable absolutely continuous. Moreover $\eE[\varepsilon_i]=0$ and $\eE[\varepsilon^2_i]=\sigma^2$.\\
 
 The assumption (H1) is standard for a change-point model, see \cite{Bai:98} or \cite{Ciuperca:11a}, while (H2) is imposed when LASSO methods are used, see for example \cite{Zou:06} or \cite{Knight:Fu:00}. Assumption (H3) is classic in a regression  model. \\
  
The matrix $\textbf{C}_{n,r}$ are non-singular for all $n$ and $r$, while the matrix $\textbf{C}_r$ can be singular. Let us denote by $\textbf{C}^0_r$ the limiting matrix for the true-change-points $l^0_r$, $r=1, \cdots, K+1$. Without loss of the generality, we assume that the regressors are centered.\\

To complete the model, we shall make the usual identifiability assumption that adjacent regressions are different: $
\ef_r \neq \ef_{r+1}$, $r=1, \cdots, K$.\\

All throughout the paper, $c$ denotes a positives generic constant. For a vector $\mathbf{v}=(v_1, \cdots, v_p)$ let us denote $| \mathbf{v} | =(|v_1|, \cdots, |v_p|)$ and $| \mathbf{v} |^c =(|v_1|^c, \cdots, |v_p|^c)$. On the other hand, $\|\mathbf{v} \|$ is its Euclidean norm. All vectors are column and $v'$ denotes the transpose of $v$. For a real $x$, $[x]$ means the largest integer not larger than $x$.  \\

After from these general notations, in every section we shall give the notations used for each method.

\section{LASSO-type estimators }

 In this section we define and study the LASSO-type estimators in a change-point model.
 
\subsection{Notations} 
Under assumption (H1), between two consecutive change-points $l_{r-1}$ and $l_r$ we consider a positive sequence $\lambda_{n,(l_{r-1},l_r)} \rightarrow \infty$ for $n \rightarrow \infty$ to control the amount of regularization applied to the estimators. Based on the LASSO-type method introduced by \cite{Knight:Fu:00} in a model without breaks, let us consider the penalized sum:
\[
T_n(K,\eth_1,\eth_2 )\equiv\sum^n_{i=1} \sum^{K+1}_{r=1} \cro{(Y_i-\eX'_i \ef_r)^2+\frac{\lambda_{n,(l_{r-1},l_r)}}{l_r-l_{r-1}} \sum^p_{j=1} | \phi_{r,j}|^\gamma } \e1_{l_{r_1} \leq i < l_r},
\]
For the tuning parameter $\lambda_{n,(l_{r-1},l_r)}=O(l_r-l_{r-1})^{1/2}$ and  $\gamma >0$, denote then the minimum of the penalized sum of the residuals squared for each fixed breaks $l_1, \cdots, l_{K+1}$:
\begin{equation}
\label{defS}
S(l_1, \cdots, l_K)\equiv\inf_{\eth_1} T_n(K,\eth_1,\eth_2 ),
\end{equation}
with $l_0=1$ and $l_{K+1}=n$. For $K$ fixed, we define the LASSO-type (or Bridge) estimator of $(\eth^0_1,\eth^0_2)$ as a point 
\[(\hat \eth^s_{1n},\hat \eth^s_{2n})=\argmin_{(\eth_1,\eth_2)} T_n(K,\eth_1,\eth_2 ).\]
 More exactly we denote: $\hat \eth^s_{1n}\equiv (\hat \ef^s_{1n}, \cdots, \hat \ef^s_{K+1,n})$ the regression parameter estimator and $\hat \eth_{2n} \equiv (\hat l_1^s, \cdots, \hat l_K^s) $ the LASSO-type estimator for the change-points. Obviously \\
 $\hat \eth_{2n} =\argmin_{(l_1, \cdots, l_K)\in \mathbb{N}^K} S(l_1, \cdots, l_K)$. \\

We remark two particular cases: for $\gamma=2$ we obtain the ridge estimator and for $\gamma=1$ the LASSO estimator. \\

The construction of the estimators has two stages: first we search the regression parameters estimators and then we localize the change-points. Then, first, between every break $l_{r-1}$ and $l_r$, we calculate the LASSO-type estimator of $\ef_r$ by: 
\[\hat \ef^s_{(l_{r-1},l_r)}=\argmin_\ef [\sum^{l_r}_{i=l_{r-1}+1}(Y_i-\eX_i' \ef)^2 +\lambda_{n;(l_{r-1},l_r)} \sum^p_{k=1} |\phi_{,k} |^\gamma],\]
 $\hat \phi^s_{(l_{r-1},l_r),k}$ its $k$th component and the corresponding forecast $\hat Y^s_{(l_{r-1},l_r),j}=\eX_j' \hat \ef^s_{(l_{r-1},l_r)}$. The regression parameters estimators are: $\hat \eth^s_{1n}(\eth_2)=(\hat \eth^s_{(l_0,l_1)}, \hat \eth^s_{{(l_1,l_2)}},\cdots, \hat \eth^s_{(l_K,l_{K+1})} )$. After we calculate the change-points estimators: 
 \[
 \hat \eth^s_{2n}=\argmin_{\eth_2 \in \N^K} T_n(K,\hat \eth^s_{1n}(\eth_2),\eth_2)
 \] 
Also $\hat \eth^s_{1n}=\eth^s_{1n}(\hat \eth^s_{2n})$.\\ 
 
  Note that in order  to take into account the sample size in every segment (phase), the tuning parameter $\lambda_{n,(l_{r-1},l_r)}$ varies from a segment to the other one with the interval length $l_r-l_{r-1}$. 
 For the true values of parameters $(\ef^0_1, \cdots, \ef^0_{K+1},$ $ l^0_1, \cdots, l^0_K)$, we consider the equivalent sum of (\ref{defS}):  $S_0\equiv  \sum^{n}_{i=1} \varepsilon^2_i +\sum^{K+1}_{r=1}\lambda_{n;(l^0_{r-1},l^0_r)} \sum^p_{j=1} | \phi^0_{r,j} |^\gamma $, where $\phi^0_{r,j}$ denotes the  $j$th component of  $\ef^0_r$.\\

For every sample, we consider the difference of the squared errors by taking some parameter  and true parameter in the model. So, for two parameters $\ef$ and $\ef^0$, we  define the function:
\[\eta_i(\ef;\ef^0)\equiv(\varepsilon_i-X'_i(\ef-\ef^0))^2-\varepsilon_i^2, \qquad i=1, \cdots, n.\]
When a LASSO-type method is considered, for the model between two consecutive change-points $j_1 <j_2$, a  penalization term is added:
\[ \eta^s_{i;(j_1,j_2)}(\ef;\ef^0)\equiv\eta_i(\ef;\ef^0) + \frac{\lambda_{n;(j_1,j_2)}}{j_2-j_1} \cro{\sum^p_{k=1} |\phi_{,k}|^\gamma-\sum^p_{k=1} |\phi^0_{,k}|^\gamma}, \quad \textrm{for } i=j_1+1, \cdots, j_2, \]
with $\phi_{,k}$ and $\phi^0_{,k}$ the $k$th component of $\ef$, respectively of $\ef^0$. We denote  $\eta^s_i(\ef;\ef^0)\equiv\eta^s_{i;(0,n)}(\ef;\ef^0)$ and $\lambda_n\equiv\lambda_{n;(0,n)}$. \\

Once the LASSO-type estimators and notations being introduced, we can study the asymptotic behavior of the estimators. First, supposing that the change-points number is known,  we consider the corresponding convergence rate, which will allow us to derive the asymptotic distribution of the estimators.  We propose a consistently estimator for the case when the  number $K$ is unknown. \\

\subsection{Asymptotic behavior}

Following result implies that the penalized sum $S$ of (\ref{defS}), optimized with respect to the regression parameters, is of order $O_{\eP}(n^\alpha)$, with $\alpha >1/2$ arbitrary. 
\begin{lemma}
\label{Lemma 1}
 Under assumptions (H2), (H3), if $\lambda_{n;(j_1,j_2)}=O(n^{1/2})$, for two points $j_1,j_2 \in \{ 1, \cdots, n\}$, $\gamma >0$, $\ef^0$  the true value of the parameter, for all $\alpha >1/2$ we have:
\[
\sup_{0 \leq j_1 <j_2 \leq n} \left| \inf_\ef \sum^{j_2}_{i=j_1+1} \eta^s_{i;(j_1,j_2)} (\ef,\ef^0) \right|
=O_{\eP}\pth{\max( \lambda_{n;(j_1,j_2)},n^\alpha)}=O_{\eP}(n^\alpha). 
\]
\end{lemma}
Convergence rate of the LASSO-type estimator is of order $n^{-1/2}$ (see Knight and Fu, 2000). We study then, by the following Lemma, the penalized sum for a LASSO-type method in a model without change-points, when the regression parameters are not in a $n^{-1/2}$-neighborhood of the true value parameter.  

\begin{lemma}
\label{Lemma 2}
 Under assumptions (H2), (H3), if $\lambda_n=o(n)$, then there exists 
$\epsilon >0$ such that: $$ \liminf_{n \rightarrow \infty} \inf_{\|\ef-\ef^0 \| \geq n^{-1/2}} n^{-1} \sum^n_{i=1} \eta^s_i(\ef,\ef^0) > \epsilon.$$ 
\end{lemma}

\noindent 
This result is useful to prove the following two Lemmas. These indicate that when the data are from two different models, the LASSO-type estimator in a model with a change-point is close to the parameter of the model where most of the data came.\\
For $u \in [ \frac{3}{4}, 1)$ in assumption (H1), $v \in (0, \frac{1}{4})$, let us consider a  constant $\delta$ such that:
\begin{equation}
\label{cond1}
\delta \in (0, u -3v).
\end{equation}
For the following lemma, the size sample of model is $n_1+n_2$.
\begin{lemma}
\label{Lemma 3}
  Under assumptions (H2), (H3), for all
$n_1, n_2 \in \N$ such that $n_1 \geq n^u$, with $3/4 \leq u  \leq 1$, $ n_2 \leq n^v$, $v < 1/4$, let be the model:
\[
\begin{array}{lll}
Y_i=\eX_i' \ef^0_1+\varepsilon_i, & &i=1, \cdots, n_1 \\
Y_i=\eX_i' \ef^0_2+\varepsilon_i, & &i=n_1+1, \cdots, n_1+n_2 \\
\end{array}
\] 
with $\ef^0_1 \neq \ef^0_2$.  
We set:
$A^s_{n_1+n_2}(\ef)\equiv\sum^{n_1}_{i=1}\eta^s_{i;(0,n_1)}(\ef;\ef^0_1) +\sum^{n_1+n_2}_{i=n_1+1}\eta^s_{i;(n_1,n_1+n_2)}(\ef;\ef^0_2)$ and $\hat \ef^s_{n_1+n_2}\equiv\argmin_\ef A^s_{n_1+n_2}(\ef)$.  Under condition (\ref{cond1}) we have:\\
(i) $\| \hat \ef^s_{n_1+n_2}-\ef^0_1\| \leq n_1^{-1/2} n_1^{\frac{v+\delta}{2u}} \leq n^{-(u-v-\delta)/2}$ .\\
(ii) $\sum^{n_1}_{i=1}\eta^s_{i;(0,n_1)}(\hat \ef^s_{n_1+n_2};\ef^0_1)=O_{\eP}(1)$.
\end{lemma}
Similarly, it have that:
\begin{lemma}
\label{Lemma 4}
Consider the model:
 \[
\begin{array}{lll}
Y_i=\eX_i' \ef^0_1+\varepsilon_i, & &i=1, \cdots, k \\
Y_i=\eX_i' \ef^0_2+\varepsilon_i, & &i=k+1, \cdots, k+n_2 \\
\end{array}
\] 
with $k \in [a n_1, n_1]$, $a \in (0,1)$. Under the same conditions as in the Lemma \ref{Lemma 3} we have:\\
(i) $\sup_{a n_1 \leq k \leq n_2} \|\hat \ef^s_{k+n_2} - \ef^0_1 \| \leq n^{-(u-v-\delta)/2} $. \\
(ii) $\sup_{a n_1 \leq k \leq n_2} \left| \sum^k_{i=1} \eta_i(\hat \ef^s_{k+n_2},\ef^0_1) \right|=O_{\eP}(1)$.
\end{lemma}

The proofs of all Lemmas are given in Appendix.\\

Suppose first that the change-point number $K$ is known. We start by studying  the convergence rate of the change-point LASSO-type estimator. The proof of Theorem \ref{Theorem 1} (given in Appendix) is split into three steps. first, using Lemma \ref{Lemma 1}, we prove that the change-point estimators are to a smaller distance $n^{1/2}$ from the true values. This implies that between two consecutive true change-points $l^0_r$, $l^0_{r+1}$, there are at most two change-point estimators. This allows to prove, using also Lemma \ref{Lemma 3}, the step 2: the change-point estimators are  at a smaller distance than $n^{1/4}$ of true values. Finally, also using Lemma \ref{Lemma 4}, the Theorem \ref{Theorem 1} is proved. 

 \begin{theorem}
 \label{Theorem 1}
  Under assumptions (H1)-(H3)  we have 
 $\hat l^s_r-l^0_r=O_{\eP}(1)$, for every $r=1, \cdots, K$.
 \end{theorem}
 
Combining Theorem  \ref{Theorem 1} and the fact that the convergence rate of LASSO-type estimators in a model without change-points is $n^{-1/2}$ allow us to have immediately the convergence rate for the regression parameter estimator:
  \begin{corollary}
  \label{Theorem 2}
   Under assumptions (H1)-(H3) we have $\|\hat \ef^s_{(\hat l^s_{r},\hat l^s_{r-1})} -\ef^0_r \|= (l^0_r-l^0_{r-1})^{-1/2}O_{\eP}(1)$, for all $r=1, \cdots , K+1$.
  \end{corollary}
  
These results imply that the LASSO-type penalization does not have influence on convergence rate, which  is the same as by LS method (see Kim and Kim, 2008): of order $(l^0_r-l^0_{r-1})^{-1/2}$ for the regression parameters and $n^{-1}$ (after a change of variable) for the change-point estimators. 
 The following theorem, whose  proof is given in Appendix, using Theorem \ref{Theorem 1}, gives the asymptotic distribution of the LASSO-type estimators for the change-points. This result implies that the asymptotic distribution is the same as for change-points estimators by LS method and it depends on the errors $(\varepsilon_i)$, the design $(X_i)$ around every true change-point and on the difference $(\ef^0_{r+1}-\ef^0_{r})$.
 
 \begin{remark}
 Since the matrix of assumption (H2) can not be  full rank, the limiting distribution of the change-point estimators can not be the \textit{argmax} of a Wiener process with shift, unlike the LS estimator (see Bai and Perron, 1998).
\end{remark} 
 
  \begin{theorem}
  \label{Theorem 3}
   Under assumptions (H1)-(H3), for each $r=1, \cdots, K$, we have that $\hat l^s_r-l^0_r$ converges in distribution for $n \rightarrow \infty$ to $\argmin_{j \in \mathbb{Z}} Z_j^{(r)}$, with $Z_0^{(r)}=0$ and\\
  - for $j=1, 2, \cdots$, $Z_j^{(r)}=\sum_{i=l^0_r+1}^{l^0_r+j} \eta_i(\ef^0_r;\ef^0_{r+1}) $;\\
  - for $j=-1, -2, \cdots$, $Z_j^{(r)}=\sum_{i=l^0_r+j}^{l^0_r}\eta_i(\ef^0_{r+1};\ef^0_{r})$.
  \end{theorem}

The asymptotic distribution of LASSO-type estimators for the regression parameters is given by the following result. For the particular case $\gamma=1$, the  LASSO estimator is  asymptotically normal for the nonzero true regression coefficients. Note also that, as \cite{Knight:Fu:00}  indicate it, for $\gamma \geq 1$, the estimator  $\hat \ef^s_{(\hat l^s_{r-1},\hat l^s_r)}$ for the nonzero coefficients can be asymptotically biased.
  \begin{theorem}
  \label{Theorem 4}
   Under assumptions (H1)-(H3), if the matrix $C^0_r=\lim_{n \rightarrow \infty } \frac{1}{l^0_r-l^0_{r-1}} \sum^{l^0_r}_{i=l_{r-1}+1}\eX_i \eX_i'$ are not singular then:
  $(\hat l^s_r-\hat l^s_{r-1})^{1/2}(\hat \ef^s_{(\hat l^s_{r-1},\hat l^s_r)}-\ef^0_r)=(l^0_r-l^0_{r-1})^{1/2}(\hat \ef^s_{(\hat l^s_{r-1},\hat l^s_r)}-\ef^0_r)(1+o_{\eP}(1)) \overset{{\cal L}} {\underset{n \rightarrow \infty}{\longrightarrow}} \argmin_{\textbf{u} \in \R^p}V_r(\textbf{u})$
 with $V_r(\textbf{u})$ defined by:
 \\
 (i) if $\gamma >1$, $V_r(\textbf{u})=-2\textbf{u}'\textbf{W}_r+\textbf{u}'\textbf{C}^0_r\textbf{u}+\lambda^0_r \sum^p_{k=1} u_k sgn(\phi^0_{r,k}) |\phi^0_{r,k} |^{\gamma -1}$, with $0 \leq \lambda^0_r=\lim_{n \rightarrow \infty} \lambda_{n;(\hat l^s_{r-1},\hat l^s_r)} / (\hat l^s_r-\hat l^s_{r-1})^{1/2}$.\\
 (ii) if $\gamma=1$, $V_r(\textbf{u})=-2\textbf{u}'\textbf{W}_r+\textbf{u}'\textbf{C}_r^0 \textbf{u}+\lambda^0_r \sum^p_{k=1} [u_k sgn(\phi^0_{r,k})\e1_{\phi^0_{r,k} \neq 0} +|u_k| \e1_{\phi^0_{r,k} =0}]$, with $0 \leq \lambda^0_r=\lim_{n \rightarrow \infty} \lambda_{n;(\hat l^s_{r-1},\hat l^s_r)} / (\hat l^s_r-\hat l^s_{r-1})^{1/2}$.\\
 (iii) if $\gamma <1$ and $\lambda_{n;(\hat l^s_{r-1},\hat l^s_r)}/(\hat l^s_r-\hat l^s_{r-1})^{\gamma/2} \overset{\eP} {\underset{n \rightarrow \infty}{\longrightarrow}} \lambda^0_r \geq 0$, $V_r(\textbf{u})=-2\textbf{u}'\textbf{W}_r+\textbf{u}'\textbf{C}^0_r\textbf{u}+\lambda^0_r\sum^p_{k=1} |u_k|^\gamma \e1_{\phi^0_k=0}$,\\
In the three cases,   the random $p$-vector $\textbf{W}_r$ is the same:  ${\cal N}(\textbf{0},\sigma^2 \textbf{C}^0_r)$. 
  \end{theorem}

 We observe that the asymptotic distributions of  LASSO-type estimators for the change-points does not depend on the tuning parameter $\lambda_{n,(l_{r-1},l_r)}$ while for the regression parameters this sequence intervenes.   This is due to the constraint  $\lambda_{n;(l_{r-1},l_r)}=O(n^{1/2})$ and to the fact that the convergence rate of the change-point estimators  is of order $n^{-1}$.\\
 
  \subsection{Choice of change-point number}

Suppose now  that we don't known a priori the change-point number which is quite often the case in practice.
Then we introduce  a criterion which is going to allow  to estimate the number $K_0$ of the change-points.
Intuitively, the true $K_0$ will be the one which minimize the function $S$ (or its $\log$) with respect to $K$. In order to take into account the model complexity, the criterion is penalized by an increasing function in $K$ and model parameter number. The penalization depends also on the sample size.\\

 For a fixed change-points number  $K$, let us denote  $\hat s_K \equiv S(\hat l^s_{1,K}, \cdots, l^s_{K,K})/n$, where the function $S$ is defined by (\ref{defS}). The criterion consistency  is established below with proof in Appendix. The demonstration idea is to prove that the change-point estimator is strictly larger or smaller than the true value with a probability converging to zero. 
 
 \begin{theorem}
 \label{Theorem 5}
 Let $\hat K_n$ be the value of  $K$ that  minimizes $B(K)=n \log \hat s_K+G(K,p_K)B_n$, with $p_K=\sum^{K+1}_{r=1} \sum^p_{j=1}\e1_{\phi^0_{r,j} \neq 0}$, function $G(K,p_K)$ increasing in $K$,  $(B_n)$ a deterministic sequence such that $B_n \rightarrow \infty $, $B_n n^{-3/4} \rightarrow 0$ and $B_n n^{-1/2} \rightarrow \infty$ for $n \rightarrow \infty$. Then $\eP[\hat K_n=K_0]\rightarrow 1$ for $n \rightarrow \infty$.
 \end{theorem}

The function $G$ is a generic factor on a penalty $B_n$.  If on every segment (phase) all variables are significant, taking $G(K,p_K)=K$, we obtain the Schwarz criterion proposed by Yao (1988). \\
 Concretely in practice, we begins by finding $K$ using  the proposed criterion. Afterthat their locations and the regression parameters on the segments defined by the  change-points are   estimated by minimizing the function $T_n$ with respect to $\eth_1$. Finally, the function $S$ given by  (\ref{defS}) is minimized with respect to $\eth_2$.

\section{Adaptive LASSO}

As mentioned in Introduction, in a model without change-points,  the LASSO-type estimator has the oracle properties for $\gamma \in (0, 1)$ only, but it is not continuous. 
On the other hand, the LASSO estimator (for $\gamma=1$) is continuous but  does not satisfy the oracle properties. To remedy these inconveniences, Zou (2006) proposed an adaptive  LASSO estimator. These properties are all the most interesting in a change-point model, since they allow to select the significant predictive variables on every phrase without realizing the hypothesis tests. 
 We  propose the adaptive LASSO estimator for a model with change-points.

\subsection{Notations}
 
In a similar way to the LASSO-type method, between two consecutive change-points $l_{r-1}$ and $l_r$, we define the regression parameters estimators by adaptive LASSO method:
\begin{equation}
\label{fs*}
 \hat \ef^{s*}_{(l_{r-1},l_r)}=\argmin_\ef [\sum^{l_r}_{i=l_{r-1}+1}(Y_i-\eX_i' \ef)^2 +\lambda_{n;(l_{r-1},l_r)} \sum^p_{k=1} \hat w_{(l_{r-1},l_r),k}|\phi_{,k} |], 
\end{equation}
 with $\hat \phi^{s*}_{(l_{r-1},l_r),k}$ its $k$th component,  $\hat \ew_{(l_{r-1},l_r)}\equiv| \hat \ef_{(l_{r-1},l_r)}|^{-g}$, $\hat w_{(l_{r-1},l_r),k}$ the $k$th component of $\hat \ew_{(l_{r-1},l_r)}$ and $\hat \ef_{(l_{r-1},l_r)}$ is LS estimator of  $\ef$ calculated between $l_{r-1}$ and $l_r$. The constant $g $ is positive and will be later specified. 
 As for the LASSO-type estimator, the tuning parameter $\lambda_{n;(l_{r-1},l_r)}$ depends of the sample size in every segment. \\
 
 Now, the penalized sum for the true parameters is: $S^*_0\equiv\sum^n_{i=1}\varepsilon^2_i+\sum^{K+1}_{r=1}\lambda_{n,(l^0_{r-1},l^0_r)}\sum^p_{k=1} \hat w_{(l^0_{r-1},l^0_r),k}|\phi^0_{r,k} |$.  Consider the sum with an adaptive penalization:
$$
S^*(l_1, \cdots, l_K) \equiv \inf_{\eth_1} \sum^{K+1}_{r=1} \cro{\sum^{l_r}_{i=l_{r-1}+1} (Y_i-\eX'_i\ef_r)^2+\lambda_{n,(l_{r-1},l_r)} \sum^p_{k=1} \hat w_{(l_{r-1},l_r),k} |\phi_{r,k}|}.
$$ 

This allows us to define the adaptive LASSO estimator for the change-points: $(\hat l^{s*}_1, \cdots, \hat l^{s*}_K)\equiv \argmin_{(l_1, \cdots, l_K)\in \mathbb{N}^K} S^*(l_1, \cdots, l_K)$ and for the regression parameters, using (\ref{fs*}): $\hat \ef^{s*}_{(\hat l^{s*}_{r-1},\hat l^{s*}_r)}$, for each $r=1, \cdots, K+1$. \\

We define, between two consecutive change-points $j_1 <j_2$, for $\ef^0$ the true value of the parameters, following penalized difference:
$$\eta^{s*}_{i;(j_1,j_2)} (\ef;\ef^0)=\eta_{i} (\ef;\ef^0)+\frac{\lambda_{n,(j_1,j_2)}}{j_2-j_1} \sum^p_{k=1} \hat w_{(j_1,j_2),k}[| \phi_{,k}|-|\phi^0_{,k} |]$$
with $\phi_{,k}$ and $\phi^0_{,k}$ the $k$th component of $\ef$, respectively of $\ef^0$.\\

In order to study the oracle properties, for each two consecutive true  change-points $l^0_{r-1}$, $l^0_r$, consider the set 
\[
{\cal A}^*_{(l^0_{r-1},l^0_r)}\equiv\{k \in \{1, \cdots,p\}; \phi^0_{r,k} \neq 0\}
\]
with the index of nonzero components of the true regression parameters. Since in practical problems, we don't known the true value of the change-points, but their estimators, we consider the similar set of index, corresponding to  the adaptive LASSO estimators $\hat l^{s*}_{r-1},\hat l^{s*}_r$ of the change-points:
\[
{\cal A}^*_{n,(\hat l^{s*}_{r-1},\hat l^{s*}_r)}\equiv\{k \in \{1, \cdots,p\}; \hat \phi^{s*}_{(\hat l^{s*}_{r-1},\hat l^{s*}_r),k} \neq 0\}\]
 will be abbreviated ${\cal A}^*_r$ for convenience. For  simplicity we denote by $\ef_{{\cal A}^*_r}$ the sub-vector of $\ef$ containing the corresponding components of  ${\cal A}^*_r$.  We also denote by $C^0_{r,kj}$ the $(k,j)$th component of matrix $\textbf{C}^0_r$.\\

\subsection{Asymptotic behavior}
By the Lemmas \ref{Lemma 5} and \ref{Lemma 6} we prove that the adaptive penalization is not a bigger order than the minimized squares sum.
\begin{lemma}
\label{Lemma 5}
Let the model $\mathbf{Y}=\eeX \ef +\ee1$, with $\mathbf{Y}$ a $n \times 1$ vector of $Y_i$ and $\eeX$ a $n \times p$ matrix.  
If $\ef^0$ is the true value of the parameter $\ef$, under assumptions (H2), (H3), for  $g \in (0, \frac{1}{4})$,  $0 \leq j_1 <j_2 \leq n$, $\lambda_{n,(j_1,j_2)}=o(n^{1/2})$,   we have: $\lambda_{n,(j_1,j_2)}\hat \ew_{(j_1,j_2)}=O_{\eP}(n^{\frac{1+g}{2}})$.
\end{lemma}

Following lemmas are needed to the proofs for the adaptive LASSO results. In fact, they study the adaptive LASSO estimator in a model without change-point and they are the equivalent of the Lemmas \ref{Lemma 1}-\ref{Lemma 3}. \\
\begin{lemma}
\label{Lemma 6}
 For two points $j_1,j_2 \in \{0,1, \cdots, n\}$, $\ef^0$ the true value of the parameters, under assumptions (H2), (H3), if $\lambda_{n,(j_1,j_2)}=o(n^{1/2})$, then, for all $g \in (0, \frac{1}{4})$ we have:
 $$
\sup_{0 \leq j_1 <j_2 \leq n} \left| \inf_\ef \sum^{j_2}_{i=j_1+1} \eta^{s*}_{i;(j_1,j_2)} (\ef;\ef^0) \right|
=O_{\eP}(n^{\frac{1+g}{2}}).$$
\end{lemma}

\begin{lemma}
\label{Lemma 7}
 Under assumptions (H2), (H3), if $\lambda_n=o(n^{1/2})$, then there exists 
$\epsilon >0$ such that: $$ \liminf_{n \rightarrow \infty} \inf_{\|\ef-\ef^0 \| \geq n^{-1/2}} n^{-1} \sum^n_{i=1} \eta^{s*}_{i;(0,n)}(\ef;\ef^0) > \epsilon .$$
\end{lemma}
\begin{lemma}
\label{Lemma 8}
For all $n_1, n_2 \in \N$ such that $n_1 \geq n^u$, with $3/4 \leq u  \leq 1$, $ n_2 \leq n^v$,  $v < 1/4$, let us consider the model:
\[
\begin{array}{lll}
Y_i=\eX_i' \ef^0_1+\varepsilon_i, & &i=1, \cdots, n_1 \\
Y_i=\eX_i' \ef^0_2+\varepsilon_i, & &i=n_1+1, \cdots, n_1+n_2 \\
\end{array}
\] 
with the assumption $\ef^0_1 \neq \ef^0_2$.  
Consider $A^{s*}_{n_1+n_2}(\ef)\equiv\sum^{n_1}_{i=1}\eta^{s*}_{i;(0,n_1)}(\ef,\ef^0_1) +\sum^{n_1+n_2}_{i=n_1+1}\eta^{s*}_{i;(n_1,n_1+n_2)}(\ef;\ef^0_2)$  and $\hat \ef^{s*}_{n_1+n_2}\equiv\argmin_\ef A^{s*}_{n_1+n_2}(\ef)$.  Under the condition (\ref{cond1}) and assumptions (H2), (H3),  we have:\\
(i) $\| \hat \ef^{s*}_{n_1+n_2}-\ef^0_1\| \leq n_1^{-1/2} n_1^{\frac{v+\delta}{2u}} \leq n^{-(u-v-\delta)/2}$ .\\
(ii) $\sum^{n_1}_{i=1}\eta^{s*}_{i;(0,n_1)}(\hat \ef^{s*}_{n_1+n_2};\ef^0_1)=O_{\eP}(1)$.
\end{lemma}
 The equivalent of the Lemma \ref{Lemma 4} is also  valid for the adaptive LASSO estimators.\\

As in Section 3, let us first suppose  that the change-point number $K$ is fixed.
By the same arguments used in the proof of Theorem \ref{Theorem 1}, using now  Lemmas \ref{Lemma 6} and \ref{Lemma 8}, we have following theorem which gives the convergence rate of the adaptive LASSO estimators for the change-points. The proof is omitted.
 \begin{theorem}
 \label{Theorem 6}
  Under assumptions (H1)-(H3), for all $g \in (0,\frac{1}{4})$,   we have 
 $\hat l^{s*}_r-l^0_r=O_{\eP}(1)$, for each $r=1, \cdots, K$. 
 \end{theorem} 

In order to work in a bounded interval, we can consider $\tau^0_r \equiv l^0_r/n$ and $\hat \tau^{s*}_r \equiv \hat l^{s*}_r/n$ its estimator. Then, Theorem \ref{Theorem 6} implies that $\hat \tau^{s*}_r$ converges in probability to $\tau^0_r$ with the convergence rate $n^{-1}$. It is the same convergence rate as by the LS method (see Bai and Perron, 1998).\\

We need that $g \in (0,\frac{1}{4})$ since the equivalent of the relation (\ref{e3}), for samples between two true change-points, is $-O_{\eP}(n^{\frac{1+g}{2}})$ which must be   $\ll O_{\eP}([n^\rho])$. This  will imply that every change-point estimator is to a distance strictly smaller than $n^\rho$ with respect to the true value. The constant $1/4$ results from the supposition that the constant $u$ of assumption (H1) is larger than $3/4$. In fact, the satisfied condition by positive constant $g$ is: $g+u <1$. In a model without change-points the constant $g$ can take any positive value, while in a change-point model it depends on the distance between the change-points.\\

Theorem \ref{Theorem 6} combined with the oracle properties of adaptive LASSO estimator in a model without change-points (see Zou, 2006) yield that the convergence rate of $\hat \ef^{s*}_{(\hat l^{s*}_{r-1},\hat l^{s*}_{r})}$ to the $\ef^0_r$ is of order $(l^0_r-l^0_{r-1})^{-1/2}$, the same for the LASSO-type estimator and by LS method (see Bai and Perron, 1998).  
\begin{remark}
By similar arguments used for the proof of Theorem \ref{Theorem 3} we can prove that  the asymptotic distribution of differences  $\hat l^{s*}_{r}-l^0_r$ is the same as for the LASSO-type estimators.\\
\end{remark}

The presence of change-points  in the model makes that the important oracle properties of the adaptive LASSO estimator for the regression parameters are not obvious.\\

The following result proves that on every segment, the adaptive LASSO estimator for the regression parameters has the oracle properties: nonzero parameters estimator on each estimated segment is asymptotically normal and zero parameters are shrunk directly to 0 with a probability converging to 1.
 
 \begin{theorem}
 \label{Theorem 7}
  Under assumptions (H1)-(H3), $g \in (0,\frac{1}{4})$, if $\lambda_{n,(l^0_{r-1},l^0_r)}(l^0_r-l^0_{r-1})^{-1/2} \rightarrow 0$ and $\lambda_{n,(l^0_{r-1},l^0_r)}(l^0_r-l^0_{r-1})^{(g-1)/2} \rightarrow \infty$ for $n \rightarrow \infty$, then:\\
    (i) $
  (\hat l^{s*}_r-\hat l^{s*}_{r-1})^{1/2}(\hat \ef^{s*}_{(\hat l^s_{r-1},\hat l^{s*}_r)}-\ef^0_r)_{{\cal A}^*_r}=(l^0_r-l^0_{r-1})^{1/2}(\hat \ef^{s*}_{(\hat l^{s*}_{r-1},\hat l^{s*}_r)}-\ef^0_r)_{{\cal A}^*_r}(1+o_{\eP}(1)) \overset{{\cal L}} {\underset{n \rightarrow \infty}{\longrightarrow}} {\cal N}(\textbf{0},\sigma^2 (\eO^0_r)^{-1})$, 
 where  for $q_r=Card\{{\cal A}^*_{(l^0_{r-1},l^0_r)} \}$, $\eO^0_r=(C^0_{r,kj})_{k,j \in {\cal A}^*_{(l^0_{r-1},l^0_r)}}$ is a $q_r \times q_r$ matrix.\\
 (ii) $\lim_{n \rightarrow \infty} \eP[{\cal A}^*_{n,(\hat l^{s*}_{r-1},\hat l^{s*}_r)}={\cal A}^*_{n,(l^0_{r-1},l^0_r)}={\cal A}^*_{(l^0_{r-1},l^0_r)}]=1$.
  \end{theorem}

Theorem \ref{Theorem 7} is proved in Appendix. The demonstration is based on the Lemma \ref{Lemma 5}, the Karush-Kuhn-Tucker (KKT) conditions and the oracle properties in a model without change-points.  Note that, for nonzero coefficients their estimators are asymptotically unbiased.\\

Let us make vary the change-point number $K$. 
 Choosing  \[(\hat l^{s*}_{1,K}, \cdots, \hat l^{s*}_{K,K}) \equiv \argmin_{(l_1, \cdots, l_K)} S^*(l_1, \cdots,  l_K)\] and $\hat s^*_K \equiv S^*(\hat l^{s*}_{1,K}, \cdots, l^{s*}_{K,K})/n$, similarly to the Theorem \ref{Theorem 5} and its proof, we  can define a consistent criterion for the change-point number: $\hat K^*_n=\argmin_K(n \log \hat s^*_K+G(K,p_K)B_n)$, with function $G$ and sequence $B_n$ as in Theorem \ref{Theorem 5}.\\
 
 \begin{remark}
 Since, on every segment,  the adaptive LASSO estimator of the regression parameter has the oracle properties, we can consider instead of $p_K$ its estimator \[\hat p_K=\sum^{K+1}_{r=1} \sum^p_{j=1}\e1_{\hat \phi^{s*}_{(\hat l^{s*}_{r_1},\hat l^{s*}_r),j} \neq 0},\] which was not possible for the LASSO estimator ($\gamma=1$) studied in Section 3. Then the  adaptive LASSO criterion is more interesting numerically. 
 \end{remark}
 \begin{remark}
 If $p$ is large compared to $n$, or more precisely for $(K+1)p \geq n$, there exists at least a segment where the LS estimator cannot be calculated. Then the adaptive LASSO estimator cannot be calculated also, and in this case, the LASSO-type method must be used. 
 \end{remark}
   
\section{Simulations}

 To illustrate the theoretical results and to compare the performances of the adaptive LASSO method with classical LS method in a change-point model we perform a simulation study. By these simulations, we show the advantages of the proposed method in terms of detection of irrelevant predictive variables. The obtained results proves that the proposed method will be  very useful for an high-sized change-point model.\\ 
 All simulations were performed using the R language. To calculate the adaptive LASSO estimations, the function \textit{lqa} of the package \textit{lqa} was used. \\
 
First, the number of phases is assumed to be known.  We consider 10 latent variables $X_1, \cdots, X_{10}$ with $X_3 \sim {\cal N}(2,1)$, $X_4 \sim {\cal N}(4,1)$, $X_5 \sim {\cal N}(1,1)$ and $X_j \sim {\cal N}(0,1)$ for $j \in \{1,2,6,7,8,9,10 \}$. The models contain two change-points (three phases) and the errors are Gaussian standardized. The true values of the regression parameters (coefficients) on the three segments are respectively: $(1,0,4,0,-3,5,6,0,-1,0)$, $(0,3,-4,-3,0,1,2,-3,0,10)$, $(1,3,4,0,0,1,0,0,0,1)$. The tuning parameter $\lambda_{n;(l_{r-1},l_r)}$ is $(l_r-l_{r-1})^\rho$ and the two change-points can vary in the interval $[1,n]$. For adaptive LASSO method, various values for the parameters $g$ and $\rho$ are considered. Recall that $g$ is the power of the adaptive penalization $\hat \ew_{(l_{r-1},l_r)}\equiv| \hat \ef_{(l_{r-1},l_r)}|^{-g}$ in relation (\ref{fs*}) and $\rho$ is the power for the tuning parameter $\lambda_{n;(l_{r-1},l_r)}=(l_r-l_{r-1})^\rho$ on each interval $[l_{r-1},l_r]$ The sample size $n$ varies from 35 to 400. The classical LS method is also considered. For each model, we generated 500 Monte-Carlo random samples of size $n$. The percentage of zero coefficients incorrectly estimated to zero(true 0) and the percentage of nonzero coefficients estimated to zero(false 0) are computed (see Tables \ref{Tabl1}-\ref{Tabl5}). Since the asymptotic distribution of the change-points estimators can not be symmetric,  in each table we also give the median of the change-point estimations.\\

\begin{table}
\vspace{0.5cm}   
 \caption{Median of change-points estimations, percentage of true 0 and of false 0 by adaptive LASSO and LS methods for $n=50$, $K=2$, $l^0_1=20$, $l^0_2=35$.}
\begin{center}
\begin{tabular}{|c|cccc|c|} \hline 
$(g,\rho)=$ & $(\frac{1}{7},\frac{13}{28})$  & $(\frac{1}{6},\frac{11}{24})$& $(\frac{1}{5}, \frac{9}{20})$ &$(\frac{9}{40},\frac{2}{5})$ & LS \\ \hline
   median of $(\hat l^{s*}_1,\hat l^{s*}_2)$ & (20,35) & (20,35)   &(20,35) & (20,35) & (20,35)  \\ \hline
   $\%$ of trues 0  & 77  &  78 & 77& 77 &   0\\ \hline
   $\%$ of   false 0 & 22 & 22 & 20 &  18   & 0 \\ \hline
\end{tabular} 
\end{center}
\label{Tabl1} 
\vspace{0.5cm} 
\end{table}

\begin{table}
\vspace{0.5cm} 
\caption{Median of change-points estimations, percentage of true 0 and of false 0 by adaptive LASSO and LS methods for $n=100$, $K=2$, $l^0_1=20$, $l^0_2=85$.}
\begin{center}
\begin{tabular}{|c|ccccc|c|} \hline 
$(g,\rho)=$ & $(\frac{1}{7},\frac{13}{28})$ & $(\frac{1}{6},\frac{11}{24})$& $(\frac{1}{5}, \frac{9}{20})$ &$(\frac{9}{40},\frac{2}{5})$  & $(\frac{2}{5},\frac{7}{20})$    & LS \\ \hline
   median of $(\hat l^{s*}_1,\hat l^{s*}_2)$ & (20,85) & (20,85)   &(20,85) & (20,85) & (20,85) & (20,85) \\ \hline
   $\%$ of trues 0  & 86 & 87 &  88 & 88 &  89  & 0\\ \hline
   $\%$ of   false 0 & 18 &  18 & 17 & 14 & 12     & 0 \\ \hline
\end{tabular} 
\end{center}
\label{Tabl2} 
\vspace{0.5cm} 
\end{table}
\begin{table}
\vspace{0.5cm} 
\caption{Median of change-points estimations, percentage of true 0 and of false 0 by adaptive LASSO and LS methods for $n=400$, $K=2$, $l^0_1=20$, $l^0_2=385$.}
\begin{center}
\begin{tabular}{|c|ccc|c|} \hline 
$(g,\rho)=$  & $(\frac{1}{6},\frac{11}{24})$ & $(\frac{1}{5},\frac{9}{20})$  &$(\frac{9}{40},\frac{2}{5})$   & LS \\ \hline
   median of $(\hat l^{s*}_1,\hat l^{s*}_2)$ & (20,385)& (20,385) & (20,385) & (20,385) \\ \hline
   $\%$ of trues 0  & 90 &  90 & 90 & 0\\ \hline
   $\%$ of   false 0 & 16 & 16& 12    & 0 \\ \hline
\end{tabular} 
\end{center}
\label{Tabl3} 
\vspace{0.5cm} 
\end{table}
\begin{table}
\vspace{0.5cm} 
\caption{Median of change-points estimations, percentage of true 0 and of false 0 by adaptive LASSO and LS methods for $n=500$, $K=2$, $l^0_1=200$, $l^0_2=400$.}
\begin{center}
\begin{tabular}{|c|ccc|c|} \hline 
$(g,\rho)=$ & $(\frac{1}{6},\frac{11}{24})$& $(\frac{1}{5},\frac{9}{20})$   &$(\frac{9}{40},\frac{2}{5})$    & LS \\ \hline
   median of $(\hat l^{s*}_1,\hat l^{s*}_2)$ & (200,400)& (200,400) & (200,400) & (200,400) \\ \hline
   $\%$ of trues 0 & 99.9 & 99.9  & 100  & 0\\ \hline
   $\%$ of   false 0 & 9.5 &  8 & 4  & 0 \\ \hline
\end{tabular} 
\end{center}
\label{Tabl4} 
\vspace{0.5cm} 
\end{table}
\begin{table}
\vspace{0.5cm} 
\caption{Median of change-points estimations, percentage of true 0 and of false 0 by adaptive LASSO and LS methods for $n=1500$, $K=2$, $l^0_1=200$, $l^0_2=400$.}
\begin{center}
\begin{tabular}{|c|ccc|c|} \hline 
$(g,\rho)=$  & $(\frac{1}{6},\frac{11}{24})$& $(\frac{1}{5},\frac{9}{20})$   &$(\frac{9}{40},\frac{2}{5})$     & LS \\ \hline
   median of $(\hat l^{s*}_1,\hat l^{s*}_2)$ & (200,400) & (200,400) & (200,400)   & (200,400)\\ \hline
   $\%$ of trues 0  & 100 &100 & 100  & 0\\ \hline
   $\%$ of   false 0 & 4.9 &  3.5 & 2.3   & 0 \\ \hline
\end{tabular} 
\end{center}
\label{Tabl5} 
\vspace{0.5cm} 
\end{table}

We obtain that, if there are segments with a small sample size, the detection percentage of the true zeros is relatively low and that of detection of the false zeros is high. The same results were obtained in the simulations of Zou (2006) for a model without change-points (for a model without change-point with a sample size equal to 60, four covariates, the largest percentage obtained by Zou to detect the zeros was 73$\%$). We observe that the detection rate of the true 0 varies more with the sample size on every segment than with the parameters $g$ or $\rho$. This rate increases slightly with $g$  and it does not depend on the location of change-points: equidistant or not. Recall that the performances of criteria proposed by Wu (2008) have varied with the change-point location for fixed $n$.  In all cases, even for a small number of observations, the median of the obtained estimations coincides with the true values of the change-points. \\
Hence, when the sample size increases, the adaptive LASSO method tends to select the true model. The penalization absence means that the LS method do not exhibit this good property in the sense that all estimations are non-zero and in order to identify the 0 parameters some supplementary hypothesis test on every phase are necessary. \\ 

In order to illustrate the model selection  criterion, we now simulate a linear model with one change-point: $Y_i=\eX_i \ef^0_1 \e1_{1 \leq i < l^0_1}+\eX_i \ef^0_2 \e1_{l^0_1 \leq i \leq n}+\varepsilon_i$, $i=1, \cdots, n$, with $n=100$, $\ef^0_1=(1,0,4,0,-3,5,6,0,-1,0)$, $\ef^0_2=(1,3,4,0,0,1,0,0,0,1)$ and $l^0_1=35$. The criterion $B(K)$ for adaptive method is computed for $B_n=n^{5/8}$,  $G(K,p_K)=K$ and the parameters for the adaptive LASSO are $g=1/5$, $\rho=9/20$. For this model 200 Monte Carlo samples of size $n$ are generated for regressor $\eX$ and error $\varepsilon$. We obtain that $\argmin_{K \in \{0,1,2,3  \}}B(K)=1$ for each Monte Carlo model replication. \\

\textit{Conclusion} By these simulations we showed advantages of the proposed (adaptive LASSO) method in terms of detection of irrelevant variables in a model with change-points and also of break number estimation. For a large enough sample size the adaptive LASSO method selects the true model, independently of change-point location. On the other hand, the constant $g$ affects slightly the detection of the true 0 parameters. The change-points are correctly estimated.

\section{Appendix} 

Here we present the proofs of the results stated in Section 3 and 4. We first give the proofs of Lemmas which are useful to prove the main results. \\

 \subsection{Proofs of Lemmas}
 The first four lemmas concern the LASSO-type estimator.

\noindent {\bf Proof of Lemma \ref{Lemma 1}}
We first show that:
\begin{equation}
\label{e1}
\sup_{0 \leq j_1 <j_2 \leq n} | \inf_\ef  \sum^{j_2}_{i=j_1+1} \eta_i(\ef;\ef^0)|=O_{\eP}(n^\alpha).
\end{equation}
Since $\eE[\varepsilon_i]=0$, $\eE[\eta_i(\ef;\ef^0)] \geq 0$ and $\eta_i(\ef^0;\ef^0)=0$, it holds that: $
0 \geq \inf_{\ef} \sum^k_{i=l} \eta_i(\ef;\ef^0) \geq \inf_\ef \sum^k_{i=l} [\eta_i(\ef;\ef^0)-\eE[\eta_i(\ef;\ef^0)]]$. 
Thus $ | \inf_{\ef} \sum^k_{i=l} \eta_i(\ef;\ef^0) | \leq \sup_\ef | \sum^k_{i=l} [\eta_i(\ef;\ef^0) - \eE[\eta_i(\ef;\ef^0)]]|  $ and relation (\ref{e1}) follows as in Lemma 3 of   Bai (1998), using (H3).  For $\eta_i^s$ we have:\\ $|\eta^s_{i;(j_1,j_2)}(\ef;\ef^0)- \eta_i(\ef;\ef^0) | \leq \frac{\lambda_{n;(j_1,j_2)}}{j_2-j_1} \left| \sum^p_{k=1} |\phi_{,k} |^\gamma -  \sum^p_{k=1} |\phi^0_{,k} |^\gamma \right|$, which implies, taking into account (\ref{e1}), that: 
$
\sup_{0 \leq j_1 <j_2 \leq n} | \inf_\ef \sum^{j_2}_{i=j_1+1} \eta^s_{i;(j_1,j_2)}(\ef;\ef^0) | \leq \sup_{0 \leq j_1 <j_2 \leq n} ( | \inf_\ef \sum^{j_2}_{i=j_1+1} \eta_i(\ef;\ef^0) | +c \lambda_{n;(j_1,j_2)}) =O_{\eP}(n^\alpha)+O(n^{1/2})$.
\hspace*{\fill}$\blacksquare$ \\

\noindent {\bf Proof of Lemma \ref{Lemma 2}} Since $\lambda_n=o(n)$, then for all $\epsilon >0$, $\ef \in  \Gamma$, there exists $n_\epsilon \in \N$ such that:
\begin{equation}
\label{etoile}
 n^{-1} \lambda_n \left|\sum^p_{k=1}|\phi_{,k}|^\gamma- \sum^p_{k=1}|\phi^0_{,k}|^\gamma \right| \leq \frac{\epsilon}{2}, \qquad n \geq n_\epsilon,
\end{equation}
uniformly in $\ef$. Under (H2) and (H3), we have that $\eE[n^{-1} \sum^n_{i=1} \eta_i(\ef;\ef^0)]>0$ for $\|\ef-\ef^0 \| \geq n^{-1/2}$ and furthermore $Var[n^{-1} \sum^n_{i=1} \eta_i(\ef;\ef^0)] \leq C$. Then there exists  $\epsilon >0$ such that, with  probability 1:
\begin{equation}
\label{lemma6} 
 \liminf_{n \rightarrow \infty} \inf_{\|\ef-\ef^0 \| \geq n^{-1/2}} n^{-1} \sum^n_{i=1} \eta_i(\ef;\ef^0) > \epsilon.
\end{equation}  
On the other hand: $\inf_{\|\ef-\ef^0 \| \geq n^{-1/2}}  n^{-1} \sum^n_{i=1} \eta^s_i(\ef;\ef^0) \geq \inf_{\|\ef-\ef^0 \| \geq n^{-1/2}}  n^{-1} \sum^n_{i=1} \eta_i(\ef;\ef^0) -n^{-1} \lambda_n\sup_{\|\ef-\ef^0 \| \geq n^{-1/2}} \left|\sum^p_{k=1}|\phi^0_{,k}|^\gamma- \sum^p_{k=1}|\phi_{,k}|^\gamma \right|$. 
The Lemma follows taking into account the relations (\ref{etoile}) and (\ref{lemma6}).
\hspace*{\fill}$\blacksquare$ \\

\noindent {\bf Proof of Lemma \ref{Lemma 3}}
(i) We show first that the assertion is  true for the LS estimator: $\hat \ef_{n_1+n_2}\equiv\argmin_\ef \sum^{n_1+n_2}_{i=1}(Y_i-\eX'_i \ef)^2$. Taking into account the assumptions (H2) and (H3), we have for all $\ef \in \Gamma$, $ 
\sum^{n_1+n_2}_{i=n_1+1} \eta_i(\ef;\ef^0_2)=O_{\eP}(n^v)$. 
By way of contradiction, suppose that for the LASSO-type estimator: $\| \hat \ef_{n_1+n_2}-\ef^0_1 \| \geq n_1^{-1/2} n_1^{\frac{v+\delta}{2 u}}$. Then, since $n_1 >n^u$, we have 
 $\sum^{n_1}_{i=1} \eta_i(\hat \ef_{n_1+n_2};\ef^0_1) \geq n_1 n_1^{-1} n_1^{\frac{v+\delta}{ u}} \geq n^{v+\delta}$. Thus, taking into account that $n_2 \leq n^v$: $
 \inf_\ef \cro{\sum^{n_1}_{i=1} \eta_i(\ef;\ef^0_1)+\sum^{n_1+n_2}_{i=n_1+1} \eta_i(\ef;\ef^0_2)} \geq O_{\eP}(n^{v+\delta})+O_{\eP}(n^v)=O_{\eP}(n^{v+\delta})$. On the other hand: $
 \inf_\ef \cro{\sum^{n_1}_{i=1} \eta_i(\ef;\ef^0_1)+\sum^{n_1+n_2}_{i=n_1+1} \eta_i(\ef;\ef^0_2)} \leq \sum^{n_1+n_2}_{i=n_1+1} (\varepsilon_i-\eX_i'(\ef^0_1-\ef^0_2))^2 =O_{\eP}(n^v)$. 
Contradiction, between the two last results. Then $\| \hat \ef_{n_1+n_2} -\ef^0_1 \| \leq n^{-(u-v-\delta)/2}$.\\
 Evidently  $\sum^{n_1+n_2}_{i=n_1+1} \eta^s_{i;(n_1,n_1+n_2)}(\ef;\ef^0_2)=O_{\eP}(n^v)+o(n_2)=O_{\eP}(n^v)$. We suppose  that: $\| \hat \ef^s_{n_1+n_2}- \ef^0_1\| \geq n_1^{-1/2} n_1^{\frac{v+\delta}{2 u}}$. Then, by Lemma \ref{Lemma 2},  $
\sum^{n_1}_{i=1} \eta^s_{i;(0,n_1)}(\hat \ef_{n_1+n_2}^s;\ef^0_1) \geq O_{\eP}(\max(\lambda_{n;(0,n_1)},n_1^{\frac{v+\delta}{u}}))$ with a  strictly positive probability. 
Thus
\begin{equation}
\label{e9}
A^s_{n_1+n_2}(\ef) \geq O_{\eP}(n^v)+O_{\eP}\pth{ \max(\lambda_{n;(0,n_1)},n_1^{\frac{v+\delta}{u}})}=O_{\eP}\pth{ \max(n^{v+\delta},\lambda_{n;(0,n_1)})}
\end{equation}
But, on the other hand,  $A^s_{n_1+n_2}(\ef) \leq \sum^{n_1+n_2}_{i=n_1+1}\eta^s_{i;(n_1,n_1+n_2)}(\ef^0_1;\ef^0_2)=O_{\eP}(n_2)=O_{\eP}(n^v)$ with probability 1, what is contradictory with (\ref{e9}).\\
(ii) 
Let us denote: 
$Z_n(\ef)\equiv\sum^{n_1}_{i=1}\eta_i(\ef;\ef^0_1)$, 
$t_n(\ef)\equiv\sum^{n_1+n_2}_{i=n_1+1} [(\varepsilon_i-\eX_i'(\ef-\ef^0_2))^2-(\varepsilon_i-\eX_i'(\ef^0_1-\ef^0_2))^2 ]$, $\hat \ef_{n_1}$ is the LS estimator of $\ef^0_{1}$ calculated for $i=1, \cdots, n_1$.
For $t_n$, we use the inequality $|a^2-b^2| \leq (a-b)^2$, assumption (H2), claim (i), condition (\ref{cond1}) and we obtain:
$|t_n(\hat \ef_{n_1+n_2})| \leq O_{\eP}(n^{- (u-v-\delta)} n^v)=O_{\eP}(n^{ - (u-2v-\delta)}) =o_{\eP}(1)$. We apply this result in the following: 
$0=Z_n(\ef^0_1)=t_n(\ef^0_1) \geq Z_n(\hat \ef_{n_1+n_2})+t_n(\hat \ef_{n_1+n_2}) \geq \inf_\ef Z_n(\ef)-|o_{\eP}(1)|$.
Thus $| Z_n(\hat \ef_{n_1+n_2})| \leq | \inf_\ef Z_n(\ef)|+o_{\eP}(1)$. Under assumptions (H2) and (H3) we have:\\ $\inf_\ef Z_n(\ef)= (\sqrt{n_1}\| \hat \ef_{n_1}-\ef^0_1\|)^2 n_1^{-1} \sum_{i=1}^{n_1} \|\eX_i\|^2 - 2 \sqrt{n_1}(\hat \ef_{n_1}-\ef^0_1)n_1^{-1/2} \sum^{n_1}_{i=1}\varepsilon_i \eX_i$\\$
=O_{\eP}(1)O(1)-O_{\eP}(1)o_{\eP}(1)=O_{\eP}(1)$ and $|Z_n(\hat \ef_{n_1+n_2})|=O_{\eP}(1)$. \\
Let us denote  $t^s_n(\ef) \equiv t_n(\ef)+\lambda_{n;(n_1,n_1+n_2)} \cro{ \sum^p_{k=1}(|\phi_{,k} |^\gamma-| \phi^0_{1,k}|^\gamma)}$, \\
$Z^s_n(\ef) \equiv Z_n(\ef)+\lambda_{n;(0,n_1)} \cro{ \sum^p_{k=1}(|\phi_{,k} |^\gamma-| \phi^0_{1,k}|^\gamma)}$. 
We can write 
 $A^s_{n_1+n_2}(\ef)=Z^s_n(\ef)+t^s_n(\ef)-\sum^{n_1+n_2}_{i=n_1+1} [\varepsilon^2_i-(\varepsilon_i-\eX'_i(\ef^0_1-\ef^0_2))^2]+\lambda_{n;(n_1,n_1+n_2)}\sum^p_{k=1}(|\phi^0_{1,k} |^\gamma-| \phi^0_{2,k}|^\gamma)$. Then $\hat \ef^s_{n_1+n_2}=\argmin_\ef (Z^s_n(\ef)+t^s_n(\ef))=\argmin_\ef A^s_{n_1+n_2}(\ef)$. But 
\[ 
| t^s_n(\hat \ef^s_{n_1+n_2})| \leq \sum^{n_1+n_2}_{i=n_1+1} \|\hat \ef^s_{n_1+n_2} - \ef^0_1 \|^2 \eX_i' \eX_i +\lambda_{n;(n_1,n_1+n_2)} \left| \sum^p_{k=1}(|\hat \phi^s_{n_1+n_2,k} |^\gamma-| \phi^0_{1,k}|^\gamma) \right|
\]
and by claim (i), under condition  (\ref{cond1}):
\[
=O_{\eP}(n^{-(u-v-\delta)} n^v)+O(n^{v/2}) O_{\eP}(\| \hat \phi^s_{n_1+n_2}- \phi^0_1\|)=O_{\eP}(n^{-\frac{u-2v-\delta}{2}})=o_{\eP}(1).
\]
Besides $ Z^s_n(\ef^0_1)=t^s_n(\ef^0_1)=0$, thus:
$0 \geq \inf_\ef (Z^s_n(\ef)+t^s_n(\ef))=Z^s_n(\hat \ef^s_{n_1+n_2})+t^s_n(\hat \ef^s_{n_1+n_2}) = Z^s_n(\hat \ef^s_{n_1+n_2})- | o_{\eP}(1)| \geq \inf_\ef Z^s_n(\ef) - | o_{\eP}(1)|$.
Then
\begin{equation}
\label{e10}
 |Z^s_n(\hat \ef^s_{n_1+n_2}) | \leq |\inf_\ef Z^s_n(\ef) | +o_{\eP}(1).
\end{equation}
On the other hand:
$\inf_\ef Z^s_n(\ef) \leq \inf_\ef Z_n(\ef) +\lambda_{n;(1,n_1)} \inf_\ef \cro{\sum^p_{k=1} \pth{|\phi_{,k} |^\gamma -| \phi^0_{1,k}|^\gamma}}$. 
But $\inf_\ef \cro{\sum^p_{k=1} \pth{|\phi_{,k} |^\gamma -| \phi^0_{1,k}|^\gamma}}  \leq \sum^p_{k=1} \pth{|\hat \phi_{n_1,k} |^\gamma -| \phi^0_{1,k}|^\gamma}= - O_{\eP}(\|\hat \ef_{n_1} - \ef^0_1 \|)=- O_{\eP} \pth{n_1^{-1/2}}$. 
Since $\inf_\ef Z_n(\ef)=O_{\eP}(1)$, we get: 
$ \inf_\ef Z^s_n(\ef) \leq -O_{\eP}(1)$.
Replacing  this last relation in (\ref{e10}) we obtain: $ | Z^s_n(\hat \ef^s_{n_1+n_2}) | =O_{\eP}(1)$. 
\hspace*{\fill}$\blacksquare$ \\

Recall that $\hat {\mathbf{w}}_{(j_1,j_2)}=|\hat \ef_{(j_1,j_2)} |^{-g}$, where $\hat \ef_{(j_1,j_2)}$ is the LS estimator of $\ef$ calculated to the samples $j_1+1, \cdots, j_2$. \\
\noindent {\bf Proof of Lemma \ref{Lemma 5}}
 Since the distribution of  $\varepsilon$ is  absolutely  continuous, then $\eP[\hat \phi_{(j_1,j_2),k}=0]=0$, for every  $k \in \{1, \cdots, p \}$. Consider, in this case  by convention $\lambda_{n,(j_1,j_2)}\hat \phi_{(j_1,j_2),k}=0$ (see P\"{o}tscher and Schneider, 2009). So only the case $\hat \phi_{(j_1,j_2),k}\neq 0$ is considered. \\
Let us denote $\eeX_{(j_1,j_2)}$ a $(j_2-j_1) \times p$  sub-matrix of $\eeX$ corresponding of the samples $i=j_1+1, \cdots, j_2$, $\eX^{k}_{(j_1,j_2)}$ denotes its $k$th column and $\mathbf{Y}_{(j_1,j_2)}\equiv(Y_i)_{j_1 \leq i \leq j_2}$. \\
For $j_1,j_2 \in \{ 0,1, \cdots , n\}$, let be the set: ${\cal A}^*_{n;(j_1,j_2)}\equiv\{k \in \{1, \cdots,p\}; \hat \phi^{s*}_{(j_1,j_2),k} \neq 0\}$ with the index of nonzero components of the adaptive LASSO estimator of $\ef$ calculated to the samples $j_1+1, \cdots, j_2$. In order to prove the Lemma, we consider two possible cases for an index: it belongs or not to this set.  \\
\underline{If $k \in {\cal A}^*_{n;(j_1,j_2)}$}, then by the Karush-Kuhn-Tucker (KKT) conditions (supposing $sgn(\hat \phi^{s*}_{(j_1,j_2),k})=+$), we have: $2^{-1} \lambda_{n,(j_1,j_2)}\hat w_{(j_1,j_2),k}={\eX^{k'}_{(j_1,j_2)}}(\eY_{(j_1,j_2)}-\eeX_{(j_1,j_2)}'\hat \ef^{s*}_{(j_1,j_2)} )={\eX^{k'}_{(j_1,j_2)}}(\eeps_{(j_1,j_2)}-\eeX'_{(j_1,j_2)}(\hat \ef^{s*}_{(j_1,j_2)}-\ef^0))=O_{\eP}(n^{1/2})- O_{\eP}(n^{-1/2}) O(n)=O_{\eP}(n^{1/2})$.\\
\underline{If $k \not \in {\cal A}^*_{n;(j_1,j_2)}$}, we have that $n^{1/2} \hat \ef_{(j_1,j_2)}$ converges  in law, for $ n \rightarrow \infty$, to some centered normal distribution. Then
$\lambda_{n,(j_1,j_2)}\hat w_{(j_1,j_2),k}= \frac{\lambda_{n,(j_1,j_2)}n^{g/2}}{\pth{ {\sqrt n} |\hat \phi_{(j_1,j_2),k} |}^g} =O_{\eP}(n^{\frac{1+g}{2}})$.
\hspace*{\fill}$\blacksquare$ \\

\noindent {\bf Proof of Lemma \ref{Lemma 6}}
Obviously $
|\eta^{s*}_{i;(j_1,j_2)} (\ef;\ef^0) -\eta_{i} (\ef;\ef^0)| \leq \frac{\lambda_{n,(j_1,j_2)}}{j_2-j_1} | \sum^p_{k=1} \hat w_{(j_1,j_2),k}[| \phi_{,k}|-|\phi^0_{,k} |]|$. 
Using Lemma \ref{Lemma 5} we obtain:
\[
\sup_{0 \leq j_1 <j_2 \leq n} \left| \inf_\ef \sum^{j_2}_{i=j_1+1} \eta^{s*}_{i;(j_1,j_2)} (\ef;\ef^0) \right| \leq O_{\eP}(n^\alpha)+O_{\eP}(n^{\frac{1+g}{2}})=O_{\eP}(n^{\frac{1+g}{2}}).
\]
\hspace*{\fill}$\blacksquare$ \\

\noindent {\bf Proof of Lemma \ref{Lemma 7}}
From Lemma \ref{Lemma 5}: $
n^{-1} \lambda_n |\sum^p_{k=1} \hat w_{(0,n),k} [\phi_{,k} -\phi^0_{,k}]| \leq \epsilon/2$, $\forall n>n_\epsilon$. 
The rest of proof is similar to that of Lemma \ref{Lemma 2}.
 \hspace*{\fill}$\blacksquare$ \\

\noindent {\bf Proof of Lemma \ref{Lemma 8}}
(i) If $\| \hat \ef^{s*}_{n_1+n_2}-\ef^0_1\| \geq n_1^{-\frac{1}{2}} n_1^{\frac{v+\delta}{2u}}$, then we have: $\sum^{n_1}_{i=1}\eta^{s*}_{i;(0,n_1)}(\ef;\ef^0_1) \geq O_{\eP} \pth{\max( n_1^{\frac{v+\delta}{u}},n_1^{\frac{1+g}{2}})}=O_{\eP} \pth{\max(n^{v+\delta},n^{u \frac{1+g}{2}}) }$ and $\sum^{n_1+n_2}_{i=n_1+1}\eta^{s*}_{i;(n_1,n_1+n_2)}(\ef;\ef^0_2)=O_{\eP}(n^v)+O_{\eP}(n_2^{\frac{1+g}{2}})=O_{\eP}(n^v)+O_{\eP}(n^{v\frac{1+g}{2}})=O_{\eP}(n^v)$. These imply that $A^{s*}_{n_1+n_2}(\ef) \geq O_{\eP}\pth{\max(n^{v+\delta},n^{u \frac{1+g}{2}}) }$ with a strictly  positive probability. On the other hand, taking into account Lemma \ref{Lemma 5}, $A^{s*}_{n_1+n_2}(\ef)\leq  \sum^{n_1+n_2}_{i=n_1+1}\eta^{s*}_{i;(n_1,n_1+n_2)}(\ef^0_1;\ef^0_2)=O_{\eP}(n_2)+O_{\eP}(n_2^{\frac{1+g}{2}})=O_{\eP}(n_2)=O_{\eP}(n^v)$. Contradiction.\\
(ii) We set:
$t^{s*}_n(\ef)\equiv t_n(\ef)+\lambda_{n;(n_1,n_1+n_2)} \sum^p_{k=1}\hat w_{(n_1,n_1+n_2),k}(|\phi_{,k} |-| \phi^0_{1,k}|)$, $Z^{s*}_n(\ef)\equiv Z_n(\ef)+\lambda_{n;(0,n_1)} \sum^p_{k=1}\hat w_{(0,n_1),k}(|\phi_{,k} |-| \phi^0_{1,k}|)$, with $Z_n(\ef)$, $t_n(\ef)$ as in the proof of Lemma \ref{Lemma 3}.  
Then $A^{s*}_{n_1+n_2}(\ef)=Z^{s*}_n(\ef)+t^{s*}_n(\ef)-\sum^{n_1+n_2}_{i=n_1+1} [\varepsilon_i^2-(\varepsilon_i-\eX'_i(\ef^0_1-\ef^0_2))^2]+\lambda_{n;(n_1,n_1+n_2)} \sum^p_{k=1}\hat w_{(n_1,n_1+n_2),k}(|\phi^0_{1,k} |-| \phi^0_{2,k}|)$.
The rest of proof  is similar to that of Lemma \ref{Lemma 3}.
 \hspace*{\fill}$\blacksquare$ \\

\subsection{Proofs of theorems}
\noindent {\bf Proof of Theorem \ref{Theorem 1}}
The proof is split  into three steps, using the same technique as in the paper Ciuperca (2011a):\\
\textit{\underline{Step 1.}} We prove that, under (H1)-(H3), with probability approaching 1, the change-point estimators are to a smaller distance than $[n^{1/2}]$ from the true values. More precisely,  for $\rho \in (\alpha,1)$ with   $\alpha >1/2$ as in the Lemma \ref{Lemma 1} we have:
\begin{equation}
\label{theorem3}
\forall r=1, \cdots, K \qquad \eP[| \hat l^s_r-l_r^0 | >[n^\rho]] \rightarrow 0, \quad n \rightarrow \infty 
\end{equation}
For this, let us study what happens if we assume that there exists  $l^0_r$ such that $|l_t-l^0_r| >[n^\rho]$, for each $ t=0,1, \cdots, K+1$. In this case, between $l^0_r-[n^\rho]$ and $l^0_r+[n^\rho]$ there is not any point $l_1, \cdots, l_K$. Then, let be the set: ${\cal L}(\rho)\equiv \{(l_1, \cdots, l_K); 0 < l_1 < \cdots <l_K < n, \sum^K_{r=1} |l_r-l^0_r| \leq [n^\rho]  \}$, with $\rho \in (\alpha, 1)$ and let us consider $(l_1, \cdots, l_K) \in {\cal L}^c_r(\rho)$, with ${\cal L}^c_r(\rho)\equiv \{(l_1, \cdots, l_K); | l_t-l^0_r| >[n^\rho], \forall t=1, \cdots, K  \}$. Thus, for all $\gamma >0$, we have:
\begin{equation}
\label{e0}
S(l_1, \cdots, l_K) \geq S(l_1,\cdots, l_K, l^0_1, \cdots , l^0_{r-1}, l^0_r-[n^\rho], l^0_r+[n^\rho], l^0_{r+1}, \cdots, l^0_K )  \equiv \sum^{K+2}_{t=1}L_t,
\end{equation}
$L_t$ will be defined later. 
On the other hand $S(\hat l^s_1, \cdots, \hat l^s_K) \leq S_0$ with probability one. Recall that $S_0\equiv  \sum^{n}_{i=1} \varepsilon^2_i +\sum^{K+1}_{r=1}\lambda_{n;(l^0_{r-1},l^0_r)} \sum^p_{j=1} | \phi^0_{r,j} |^\gamma $ and the definition of $S$ is given by (\ref{defS}). 
For all   $t \in \{1, \cdots, r-1, r+1, \cdots, K+2 \}$ let us consider the points $k_{1,t} < \cdots < k_{J(t),t}\equiv \{l_1, \cdots, l_K \} \cap \{ j; l^0_{t-1} < j \leq l^0_t \}$ and define:
\[  
L_t \equiv \sum^{J(t)+1}_{j=1} \min_{\ef_j} \pth{ \sum^{k_{j,t}}_{i=k_{j-1,t}+1}(\varepsilon_i -\eX_i'(\ef_j-\ef^0_t))^2+ \lambda_{n;(k_{j-1,t},k_{j,t})} \sum^p_{k=1} | \phi_{j,k}|^\gamma}.
\] 
Hence, due to the fact the $\lambda_{n;(k_{j-1,t},k_{j,t})} =o(k_{j,t}-k_{j-1,t})$ and using Lemma \ref{Lemma 1}:
\begin{equation}
\label{e3}
\begin{array}{l}
0 \geq L_t- \sum^{J(t)+1}_{j=1}  \pth{ \sum^{k_{j,t}}_{i=k_{j-1,t}+1} \varepsilon^2_i +\lambda_{n;(k_{j-1,t},k_{j,t})} \sum^p_{k=1} | \phi^0_{j,k}|^\gamma} \\
 \geq -2(K+1) \sup_{1 \leq l < j \leq n}  \left| \inf_\phi \cro{ \sum^j_{i=l+1} \eta^s_{i;(l,j)}(\ef;\ef^0) }\right|   =-O_{\eP}(n^\alpha),
\end{array}
\end{equation}
with $\ef^0$ of relation (\ref{e3}) one of the true parameters $\ef^0_r$, $r=1, \cdots, K+1$. For the samples between  $l^0_r-[n^\rho]$ and $l^0_r+[n^\rho]$ we have:
\begin{equation}
\label{e4} 
\begin{array}{l}
L_r-\sum^{l^0_r+[n^\rho]}_{i=l^0_r-[n^\rho]+1} \varepsilon^2_i - \sum^{J(r)+1}_{j=1}\lambda_{n;(k_{j-1,r},k_{j,r})} \sum^p_{k=1} | \phi^0_{r,k}|^\gamma - \sum^{J(r)}_{j=1}\lambda_{n;(k_{j-1,r},k_{j,r})} \sum^p_{k=1} | \phi^0_{r+1,k}|^\gamma \\
=\min_\phi \left\{ \sum^{l^0_r}_{i=l^0_r-[n^\rho]+1} \eta^s_{i;(l^0_r-[n^\rho],l^0_r)}(\ef;\ef^0_r)+ \sum^{l^0_r+[n^\rho]}_{i=l^0_r+1} \eta^s_{i;(l^0_r,l^0_r+[n^\rho])} (\ef;\ef^0_{r+1}) \right\}.
\end{array}
\end{equation}
Since to left and to right of each change-point the models are different, we suppose $\|\ef -\ef^0_r \| >c$. 
But
 $\sum^{l^0_r}_{i=l^0_r-[n^\rho]+1} \eta_i(\ef;\ef^0_r) =O_{\eP}([n^\rho])\gg O(n^{1/2})$. Since $\lambda_{n;(l^0_r-[n^\rho],l^0_r)}=O(n^{1/2})$ and taking into account the relation between  $\eta_i$ and $\eta^s_i$, we obtain that:  (\ref{e4}) is of order $O_{\eP}([n^\rho]) >0$.
Then, for (\ref{e0}), using (\ref{e3}) and the last relation, we obtain: $S(l_1, \cdots, l_K) \geq -O_{\eP}(n^\alpha) +O_{\eP}([n^\rho])+S_0$.
This last relation implies:\\ $\eP [\min_{(l_1, \cdots, l_K) \in {\cal L}^c_r(\rho)} S(l_1, \cdots, l_K) > S_0 ] \rightarrow 1$ for $n \rightarrow \infty$, and relation (\ref{theorem3}) follows.\\
\textit{\underline{Step 2.}} We prove now that the change-point estimators are at a smaller distance than $n^{1/4}$ of true values: for all $\nu <1/4$ we have: $\eP[|\hat l^s_r-l^0_r | >n^\nu] \rightarrow 0$, for $n \rightarrow \infty$, for every $r=1, \cdots, K$.\\
Therefore the step 1, the change-point estimators belong in the set ${\cal L}(\rho)$,  with probability tending to  1 for $n \rightarrow \infty$. For a $r \in \{1, \cdots, K \}$ let be the subset of ${\cal L}(\rho)$:
 $$
 {\cal L}^c_{r}(\nu) \equiv \{ (l_1, \cdots, l_K) \in {\cal L}(\rho); |l_t-l^0_r|   > n^\nu, t=1, \cdots, K  \}.
 $$
 For $(l_1, \cdots, l_k) \in {\cal L}^c_{r}(\nu)$, we have that:
 $ S(l_1, \cdots, l_k) \geq S(l_1, \cdots, l_K, l^0_1, \cdots, l^0_{r-1}, l^0_r-[n^\nu], l^0_r+[n^\nu], l^0_{r+1}, \cdots, l^0_K)$.
 For $t \neq r-1,r$, by assumption (H1), using the step 1, we have that there are at most two points $l_t$ and $l_{t+1}$ between $l^0_t$ and $l^0_{t+1}$. Suppose that there are two points $l_t$ and $l_{t+1}$ between $l^0_t$ and $l^0_{t+1}$. If there is a single  point or no point the approach is the same.  
 \[
 \begin{array}{lll}
 D(l^0_t,l^0_{t+1}) & =& \inf_\ef \left\{\sum^{l_t}_{i=l^0_t+1}(\varepsilon_i-\eX_i'(\ef-\ef^0_{t+1}))^2 +\lambda_{n;(l^0_t,l_t)} \sum^p_{k=1} |\phi_{,k} |^\gamma \right\}\\
  & &+ \inf_\ef \acc{ \sum^{l_{t+1}}_{i=l_t+1} (\varepsilon_i-\eX_i'(\ef-\ef^0_{t+1}))^2+\lambda_{n;(l_t,l_{t+1})} \sum^p_{k=1} |\phi_{,k} |^\gamma }\\
  & &+ \inf_\ef \acc{ \sum^{l^0_{t+1}}_{i=l_t+1} (\varepsilon_i-\eX_i'(\ef-\ef^0_{t+1}))^2+\lambda_{n;(l_{t+1},l^0_{t+1})} \sum^p_{k=1} |\phi_{,k} |^\gamma }. \\
 \end{array}
 \]
 Consequently, $S(l_1, \cdots, l_K, l^0_1, \cdots, l^0_{r-1}, l^0_r-[n^\nu], l^0_r+[n^\nu], l^0_{r+1}, \cdots, l^0_K)-S_0$ can be written as: $ \sum_{t \neq r-1,r} \acc{  D(l^0_t,l^0_{t+1})-\sum^{l^0_{t+1}}_{i=l^0_t+1}\varepsilon^2_i- \lambda_{n;(l^0_t,l^0_{t+1})}  \sum^p_{k=1} |\phi^0_{t,k} |^\gamma} \\
 + \pth{ D(l^0_{r-1},l^0_{r}-[n^\nu])-\sum^{l^0_{r}-[n^\nu]}_{i=l^0_{r-1}+1}\varepsilon^2_i- \lambda_{n;(l^0_{r-1},l^0_{r}-[n^\nu])}  \sum^p_{k=1} |\phi^0_{r,k} |^\gamma} \\
  + \pth{ D(l^0_{r}+[n^\nu],l^0_{r+1})-\sum^{l^0_{r+1}}_{i=l^0_{r}+[n^\nu]+1}\varepsilon^2_i- \lambda_{n;(l^0_{r}+[n^\nu],l^0_{r+1})}  \sum^p_{k=1} |\phi^0_{r+1,k} |^\gamma }\\
   + \left(  D(l^0_{r}-[n^\nu],l^0_{r}+[n^\nu])-\sum^{l^0_{r}+[n^\nu]}_{i=l^0_{r}-[n^\nu]+1}\varepsilon^2_i- \lambda_{n;(l^0_{r}-[n^\nu],l^0_{r})}  \sum^p_{k=1} |\phi^0_{r,k} |^\gamma \right. 
  - \left. \lambda_{n;(l^0_{r},l^0_{r}+[n^\nu])}  \sum^p_{k=1} |\phi^0_{r+1,k} |^\gamma \right)\\
  \equiv A+B+C+D$,
 with $D(l^0_{r-1},l^0_{r}-[n^\nu])$, $D(l^0_{r}+[n^\nu],l^0_{r+1})$, $D(l^0_{r}-[n^\nu],l^0_{r}+[n^\nu])$ sums with the same form that  $ D(l^0_s,l^0_{s+1})$.\\
 \underline{For A}: since $|l^0_t-\hat l^s_t| \leq O_{\eP}([n^\rho])$, with $\rho < 3/4$ and under assumption (H1), taking into account Lemma \ref{Lemma 3}(ii), we have:
 \[ 
  D(l^0_t,l^0_{t+1}) - \sum^{l^0_{t+1}}_{i=l^0_t+1}\varepsilon^2_i-\lambda_{n;(l^0_t, l^0_{t+1})} \sum^p_{k=1} |\phi^0_{t,k}|^\gamma =\inf_\phi \acc{ \sum^{l_{t+1}}_{i=l_t+1} \eta^s_{i,(l_t,l_{t+1})}(\ef,\ef^0_{t+1}) }(1+o_{\eP}(1))=O_{\eP}(1).
 \] 
Similarly $B$, $C$ $=O_{\eP}(1)$.\\
 \underline{For D}: $D(l^0_{r}-[n^\nu],l^0_{r}+[n^\nu])$ is equal to 
 $
 \inf_\ef \left\{ \sum^{l^0_r}_{i=l^0_r-[n^\nu]}(\varepsilon_i-\eX_i'(\ef-\ef^0_r))^2  \right.$ \\
 $\left.+\lambda_{n;(l^0_{r}-[n^\nu],l^0_{r})}  \sum^p_{k=1} |\phi_{,k} |^\gamma
+  \sum^{l^0_r+[n^\nu]}_{i=l^0_r+1}(\varepsilon_i-\eX_i'(\ef-\ef^0_{r+1}))^2+\lambda_{n;(l^0_{r},l^0_{r}+[n^\nu])}\sum^p_{k=1} |\phi_{,k} |^\gamma  \right\}$. 
 Since $\ef^0_r \neq \ef^0_{r+1}$, if $\tilde \ef$ is the minimizer of the last relation, there exists a constant $c>0$ such that  at least one of two $\| \tilde \ef- \ef^0_r\| $ or $\| \tilde \ef- \ef^0_{r+1}\| $ is  greater than $c$. Let us suppose that it is the first one. Then, for $\|\ef-\ef^0_r \| >c$, using Lemma 1 of Babu(1989) we have: $ \sum^{l^0_r}_{i=l^0_r-[n^\nu]} \eta^s_{i;(l^0_{r}-[n^\nu],l^0_{r})}   (\ef;\ef^0_r)=\| \ef-\ef^0_r\|^2 \sum^{l^0_r}_{i=l^0_r-[n^\nu]} \| \eX_i\|^2 -2(\ef-\ef^0_r)'\sum^{l^0_r}_{i=l^0_r-[n^\nu]} \varepsilon_i \eX_i  +\lambda_{n;(l^0_{r}-[n^\nu],l^0_{r})}  \pth{  \sum^p_{k=1} |\phi_{,k} |^\gamma - \sum^p_{k=1} |\phi^0_{r,k} |^\gamma }=O_{\eP}(n^\nu)+o(n^\nu)=O_{\eP}(n^\nu)$, 
 uniformly in $\ef$. 
 So $D=O_{\eP}(n^\nu)$ and then 
\begin{equation}
\label{iiq}
\inf_{(l_1, \cdots, l_K) \in {\cal L}^c_{r}(\nu)} [S(l_1, \cdots, l_K, l^0_1, \cdots, l^0_{r-1}, l^0_r-[n^\nu], l^0_r+[n^\nu], l^0_{r+1}, \cdots, l^0_K)-S_0]>O_{\eP}([n^\nu]),
\end{equation}  
with probability tending to one  for $n \rightarrow \infty$. This involves:\\
  $  \inf_{(l_1, \cdots, l_K) \in {\cal L}^c_{r}(\nu)} S(l_1, \cdots, l_K, l^0_1, \cdots, l^0_{r-1}, l^0_r-[n^\nu], l^0_r+[n^\nu], l^0_{r+1}, \cdots, l^0_K) > S_0$   and the step  2 is proved.\\
 \textit{\underline{Step 3.}} Now, we can show the theorem: $\hat l^s_r-l^0_r=O_{\eP}(1)$ for any $r=1, \cdots, K$.\\
 Let be the set: ${\cal L}(\nu) \equiv \acc{(l_1, \cdots, l_K); |l_t-l^0_t| <[n^\nu], \forall t=1, \cdots, K}$ with $\nu <1/4$. For a ${\cal M}_1 >0$ to determine later, let also the set: 
  ${\cal L}_r(\nu, {\cal M}_1) \equiv \acc{(l_1, \cdots l_K) \in {\cal L}(\nu); l_r-l^0_r <-{\cal M}_1}$. \\
  Consider two vectors of change-points  $(m_1, \cdots, m_K) \in {\cal L}(\nu)$ and $(l_1, \cdots, l_K) \in {\cal L}_r(\nu,{\cal M}_1)$ such that  $m_t=l_t$ for $t \neq r$ and $m_r=l^0_r$.  Using the notations specified in Section 3, we can write:
  $ S(l_1, \cdots, l_K)- S(m_1, \cdots, m_K) =\{ \sum^{l_r}_{j=l_{r-1}+1} [(Y_j-\hat Y^s_{(l_{r-1},l_r),j})^2-(Y_j-\hat Y^s_{(l_{r-1},l^0_r),j})^2] 
   + \lambda_{n;(l_{r-1},l_r)} \sum^p_{k=1} |\hat \phi^s_{(l_{r-1},l_r),k} |^\gamma  -\lambda_{n;(l_{r-1},l^0_r)} \sum^p_{k=1} |\hat \phi^s_{(l_{r-1},l^0_r),k} |^\gamma \}
   +\{\sum^{l^0_r}_{j=l_r+1} [(Y_j-\hat Y^s_{(l_{r},l_{r+1}),j})^2-(Y_j-\hat Y^s_{(l_{r-1},l^0_r),j})^2] \}
  + \{ \sum^{l_{r+1}}_{j=l^0_{r+1}} [(Y_j-\hat Y^s_{(l_{r},l_{r+1}),j})^2 - (Y_j-\hat Y^s_{(l^0_{r},l_{r+1}),j})^2] 
  +  \lambda_{n;(l_r,l_{r+1})}  \sum^p_{k=1} |\hat \phi^s_{(l_r,l_{r+1}),k}|^\gamma - \lambda_{n;(l^0_r,l_{r+1})}  \sum^p_{k=1} |\hat \phi^s_{(l^0_r,l_{r+1}),k}|^\gamma \} \equiv \{I_{11}+I_{12} \} +\{ I_{21}\}+\{I_{31}+I_{32}\}$. 
  By Lemma \ref{Lemma 4}(ii) we have for $I_{11}$:
  $I_{11}=  \sum^{l_r}_{j=l_{r-1}+1} \cro{\eta_j(\hat \ef^s_{(l_{r-1},l_r)};\ef^0_r) -\eta_j(\hat \ef^s_{(l_{r-1},l^0_r)};\ef^0_r) }=O_{\eP}(1)$. By Lemma \ref{Lemma 4}(i), taking $n_1=l_r-l_{r-1} \geq [n^u]$, $u \geq 3/4$, $n_2=l^0_r-l_r$, since 
  $\lambda_{n;(l_{r-1},l_r)} =O((l_r-l_{r-1})^{1/2})$ it holds that: 
  \[
  I_{12}=\lambda_{n;(l_{r-1},l_r)} [\sum^p_{k=1} |\hat \phi^s_{(l_{r-1},l_r),k} |^\gamma -\sum^p_{k=1} |\hat \phi^s_{(l_{r-1},l^0_r),k} |^\gamma ]+ [\lambda_{n;(l_{r-1},l_r)}  - \lambda_{n;(l_{r-1},l^0_r)} ]\sum^p_{k=1} |\hat \phi^s_{(l_{r-1},l^0_r),k} |^\gamma
  \]
 \[=\lambda_{n;(l_{r-1},l_r)}O_{\eP}(\|\hat \ef^s_{(l_{r-1},l_r)} -\hat \ef^s_{(l_{r-1},l^0_r)}\|)+O_{\eP}(\lambda_{n;(l_{r-1};l_r)}- \lambda_{n;(l_{r-1};l^0_r)})=o_{\eP}(1)+o_{\eP}(l^0_r-l_r)\] $=o_{\eP}(l^0_r-l_r)$. 
 Similarly, it can be shown that $I_{31}=O_{\eP}(1)$ and $I_{32}=o_{\eP}(l^0_r-l_r)$. For $I_{21}$:
  \[ 
  I_{21}=\sum^{l^0_r}_{i=l_r+1}  [(\varepsilon_j-\eX'_i(\ef^0_{r+1}-\ef^0_r))^2-\varepsilon^2_i ]+\sum^{l^0_r}_{j=l_r+1}  [(\varepsilon_j-\eX'_i(\hat \ef^s_{(l_r,l_{r+1})}-\ef^0_r))^2
  \]
  $$
  -(\varepsilon_i-\eX'_j(\ef^0_{r+1}-\ef^0_r))^2 ]
  -\sum^{l^0_r}_{i=l_r+1} \cro{(\varepsilon_i-\eX'_i(\hat \ef^s_{(l_{r-1},l^0_r)}-\ef^0_r))^2 -\varepsilon^2_i} \equiv J_1+J_2+J_3.
  $$
  For $J_2$, Lemma  \ref{Lemma 4}(i) combined with the step  2, condition (\ref{cond1}), yield that:
  \[
  J_2=[\| \hat \ef^s_{(l_r,l_{r+1})}-\ef^0_r\|^2-\| \ef^0_{r+1}-\ef^0_r \|^2]\sum^{l^0_r}_{i=l_r+1} \| \eX_i\|^2 - 2 \pth{\hat \ef^s_{(l_r,l_{r+1})}-\ef^0_{r+1}) }'\sum^{l^0_r}_{i=l_r+1} \varepsilon_i \eX_i
  \]
    \[
  =O_{\eP}(n^{-(u-\nu-\delta)} (l^0_r-l_r))+O_{\eP}(n^{- \frac{u-\nu-\delta}{2}}(l^0_r-l_r))=O_{\eP}(n^{- \frac{u -\delta -3 \nu}{2}})=o_{\eP}(1).
  \]
  Similarly $|J_3|=o_{\eP}(1)$. 
  For $J_1$, since $\ef^0_r \neq \ef^0_{r+1}$, there exists $C_1 >0$ such that $|J_1| >C_1(l^0_r-l_r)$. 
  We choose  ${\cal M}_1 >0$ such that $|J_1|$ (then $I_{21}$ also) is bigger than $ \max(I_{12}, I_{11}, I_{31},I_{32})$. Then $\lim_{n \rightarrow \infty} \eP[(\hat l^s_1, \cdots, \hat l^s_K) \in {\cal L}_r(\nu,{\cal M}_1)]=0$ and theorem is established.
 \hspace*{\fill}$\blacksquare$ \\
  
\noindent {\bf Proof of Theorem \ref{Theorem 3}}
  Since $\eE[Z_j^{(r)}]>c j$, then for all $\tilde \delta>0$ there exists ${\cal M}_{2, \tilde \delta} >0$ such that: $\eP[|\argmin_{j \in \mathbb{Z}} Z_j^{(r)}| \leq  {\cal M}_{2,\tilde \delta}] >1 -\tilde \delta$. 
  By Theorem \ref{Theorem 1}, for all $\tilde \delta>0$, there exists ${\cal M}_{1,\tilde \delta}>0$ such that for each $r=1, \cdots, K$: $\eP[|\hat l^s_r-l^0_r | \leq {\cal M}_{1, \tilde \delta}] >1-\tilde \delta$. 
  Consider then  ${\cal M}=\max\{ {\cal M}_{1,\tilde \delta},{\cal M}_{2,\tilde \delta} \}$. For each  $|i_r| \leq {\cal M}$, let us consider: \\
  $S(l^0_1+i_0, \cdots, l^0_K+i_K)-S(l^0_1, \cdots, l^0_K)$
\[=\sum^K_{r=1} \sum^{l^0_r}_{j=l^0_{r-1}+i_{r-1}+1} \cro{\pth{Y_j -\eX_j' \hat \ef^s_{(l^0_{r-1}+i_{r-1},l^0_r+i_r)}}^2 -\pth{Y_j-\eX_j' \hat \ef^s_{(l^0_{r-1},l^0_r)}}^2}\]
 \[
  +\sum^K_{r=1} \sum^{l^0_r+i_r}_{j=l^0_{r}+1} \acc{\cro{(Y_j-\eX_j'\hat \ef^s_{(l^0_{r-1}+i_{r-1},l^0_r+i_r)})^2-\varepsilon^2_j} -\cro{(Y_j-\eX_j' \hat \ef^s_{(l^0_r,l^0_{r+1})})^2- \varepsilon^2_j} } 
  \]
  \[ +\sum^K_{r=1} \cro{\lambda_{n,(l^0_{r-1}+i_{r-1},l^0_r+i_r)}\sum^p_{k=1} | \hat \phi^s_{(l^0_{r-1}+i_{r-1},l^0_r+i_r),k}|^\gamma -\lambda_{n,(l^0_{r-1},l^0_r)}\sum^p_{k=1} | \hat \phi^s_{(l^0_{r-1},l^0_r),k}|^\gamma} \]
 and by Corollary \ref{Theorem 2}: $=o_{\eP}(1)+\sum^K_{r=1}Z_j^{(r)}(1+o_{\eP}(1))+o_{\eP}(1)$.  
 So $S(l^0_1+i_0, \cdots, l^0_K+i_K)-S(l^0_1, \cdots, l^0_K)$ converges jointly in $i_r$ in distribution to $\sum^K_{r=1}Z_j^{(r)}$ for $n \rightarrow \infty$. Since we have independent set of random variables around each change-point, the theorem follows.
 \hspace*{\fill}$\blacksquare$ \\
  
\noindent {\bf Proof  of Theorem \ref{Theorem 4}}
  We combine  Theorems 2 and 3 of Knight and Fu (2000) with Corollary \ref{Theorem 2} and Theorem \ref{Theorem 3} of this paper. 
 \hspace*{\fill}$\blacksquare$ \\
  
\noindent {\bf Proof of Theorem \ref{Theorem 5}}
Consider first $K <K_0$. By Lemma  \ref{Lemma 1} we have for $\alpha \in (1/2; 3/4)$: $0 > S(\hat l^s_{1,K_0}, \cdots, \hat l^s_{K_0,K_0}) -S_0=O_{\eP}(n^\alpha)$. So $\hat s_{K_0}=n^{-1} S(\hat l^s_{1,K_0}, \cdots, \hat l^s_{K_0,K_0})- n^{-1} S_0+n^{-1} S_0=O_{\eP}(n^{\alpha-1})+n^{-1} \sum^n_{i=1} \varepsilon^2_i +o_{\eP}(1)\overset{{\eP}} {\underset{n \rightarrow \infty}{\longrightarrow}} \eE[\sigma^2] $. On the other hand, since the distance between two consecutive  change-points is at least  $n^{3/4}$, we prove in the same way as for relation (\ref{iiq}), that: $S(\hat l^s_{1,K}, \cdots,\hat l^s_{K,K}) -S_0 > C n^{3/4}$. Then:
$n (\hat s_K-\hat s_{K_0})\hat s^{-1}_{K_0}=$ \\$\hat s^{-1}_{K_0} \cro{(S(\hat l^s_{1,K}, \cdots,\hat l^s_{K,K}) -S_0 )-(S(\hat l^s_{1,K_0}, \cdots, \hat l^s_{K_0,K_0})-S_0)} >O_{\eP}(n^{3/4})-O_{\eP}(n^\alpha)=O_{\eP}(n^{3/4})$.
Due to fact  $B_n \ll n^{3/4}$ we have for $K <K_0$, $B(K)-B(K_0) \overset{{\eP}} {\underset{n \rightarrow \infty}{\longrightarrow}} \infty$.\\
Let us consider now that $K >K_0$. Thus: $S_0 \geq  S(\hat l^s_{1,K_0}, \cdots,\hat l^s_{K_0,K_0}) \geq S(\hat l^s_{1,K}, \cdots, \hat l^s_{K,K})$ $ \geq S(\hat l^s_{1,K}, \cdots,\hat l^s_{K,K}, l_1^0,\cdots, l^0_{K_0}) \geq S_0-O_{\eP}(n^\nu)$, with $\nu <1/4$. The last inequality is obtained by similar calculations of the  Theorem \ref{Theorem 1} proof, step 2. Then $0 \leq \hat s_{K_0}- \hat s_K=O_{\eP}(n^{\nu -1})$, which implies immediately  that: $n \log \hat s_{K_0}-n \log \hat s_K=O_{\eP}(n^\nu)$. Taking into account that the function $G$ is increasing in $K$ and that  $n^\alpha  \ll B_n \ll n^{3/4}$, we have for the criterion $B(K)-B(K_0)=-O_{\eP}(n^\nu)+G(K,p)B_n-G(K_0,p_0)B_n > O_{\eP}(n^\alpha)$. Then, for $n \rightarrow \infty$, $\eP[\hat K_n >K_0] \rightarrow 0$.
 \hspace*{\fill}$\blacksquare$ \\

\noindent {\bf Proof of Theorem \ref{Theorem 7}}
  Observe that the claim (i) follows immediately by Theorem 2 of Zou (2006) and Theorem \ref{Theorem 6}.\\
  (ii) Recall first the definition of the two sets: ${\cal A}^*_{(l^0_{r-1},l^0_r)}\equiv\{k \in \{1, \cdots,p\}; \phi^0_{r,k} \neq 0\}$ and  ${\cal A}^*_{n,(\hat l^{s*}_{r-1},\hat l^{s*}_r)}\equiv\{k \in \{1, \cdots,p\}; \hat \phi^{s*}_{(\hat l^{s*}_{r-1},\hat l^{s*}_r),k} \neq 0\}$.  By Theorem 2 of Zou (2006) the adaptive LASSO estimator in a model without change-points has the oracle properties: $\lim_{n \rightarrow \infty} \eP[{\cal A}^*_{n,(l^0_{r-1},l^0_r)}={\cal A}^*_{(l^0_{r-1},l^0_r)}]=1$. It remains to prove that: $\lim_{n \rightarrow \infty} \eP[{\cal A}^*_{n,(\hat l^{s*}_{r-1},\hat l^{s*}_r)}={\cal A}^*_{(l^0_{r-1},l^0_r)}]=1$. The asymptotic normality of estimators implies that: $\phi^0_{r,k}- \hat \phi^{s*}_{(\hat l^{s*}_{r-1},\hat l^{s*}_r),k} \overset{\eP} {\underset{n \rightarrow \infty}{\longrightarrow}}  0$, $\forall k \in {\cal A}^*_{n,(\hat l^{s*}_{r-1},\hat l^{s*}_r)}$. Then $\lim_{n \rightarrow \infty} \eP[{\cal A}^*_{n,(\hat l^{s*}_{r-1},\hat l^{s*}_r)} \subseteq {\cal A}^*_{(l^0_{r-1},l^0_r)}]=1$. We now prove: $\eP[\exists k \in\{1, \cdots, p \}; k \not \in {\cal A}^*_{(l^0_{r-1},l^0_r)} \textrm{ and } k \in  {\cal A}^*_{n,(\hat l^{s*}_{r-1},\hat l^{s*}_r)} ]\rightarrow  0$, for $n \rightarrow \infty$. Using the notations given in the proof of Lemma  \ref{Lemma 5}, by KKT optimality conditions: 
  $  2 \lambda_{n, (\hat l^{s*}_{r-1},\hat l^{s*}_r)}  \hat w_{(\hat l^{s*}_{r-1},\hat l^{s*}_r),k}=\eX^{k'}_{(\hat l^{s*}_{r-1},\hat l^{s*}_r)} (\eY_{(\hat l^{s*}_{r-1},\hat l^{s*}_r)}- \eeX'_{(\hat l^{s*}_{r-1},\hat l^{s*}_r)} \hat \ef^{s*}_{(\hat l^{s*}_{r-1},\hat l^{s*}_r)})$.
Consequence of (i): $(\hat l^{s*}_r-\hat l^{s*}_{r-1})^{1/2} \hat \phi^{s*}_{(\hat l^{s*}_{r-1},\hat l^{s*}_r),k}=O_{\eP}(1)$. Taking into account conditions imposed to $\lambda_{n, (\hat l^{s*}_{r-1},\hat l^{s*}_r)}$:  
  \[
  \frac{\lambda_{n, (\hat l^{s*}_{r-1},\hat l^{s*}_r)} \hat w_{(\hat l^{s*}_{r-1},\hat l^{s*}_r),k}}{(\hat l^{s*}_r-\hat l^{s*}_{r-1})^{1/2}}=\frac{\lambda_{n, (\hat l^{s*}_{r-1},\hat l^{s*}_r)} }{(\hat l^{s*}_r-\hat l^{s*}_{r-1})^{1/2}} (\hat l^{s*}_r-\hat l^{s*}_{r-1})^{g/2} \frac{1}{ |(\hat l^{s*}_r-\hat l^{s*}_{r-1})^{1/2} \hat \phi^{s*}_{(\hat l^{s*}_{r-1},\hat l^{s*}_r),k} |^{g}}\overset{\eP} {\underset{n \rightarrow \infty}{\longrightarrow}} \infty.
  \]
  On the other hand, in the case $(\hat l^{s*}_{r-1},\hat l^{s*}_r) \subseteq (l^0_{r-1},l^0_r)$, we have
  \[
  \frac{\eX^{k'}_{(\hat l^{s*}_{r-1},\hat l^{s*}_r)} (\eY_{(\hat l^{s*}_{r-1},\hat l^{s*}_r)}- \eeX'_{(\hat l^{s*}_{r-1},\hat l^{s*}_r)} \hat \ef^{s*}_{(\hat l^{s*}_{r-1},\hat l^{s*}_r)})}{(\hat l^{s*}_r-\hat l^{s*}_{r-1})^{1/2}}=\frac{\eX^{k'}_{(\hat l^{s*}_{r-1},\hat l^{s*}_r)} \eeX_{(\hat l^{s*}_{r-1},\hat l^{s*}_r)} (\hat l^{s*}_r-\hat l^{s*}_{r-1})^{1/2}(\ef^0_r-\hat \ef_{(\hat l^{s*}_{r-1},\hat l^{s*}_r)}) }{\hat l^{s*}_r-\hat l^{s*}_{r-1}}
  \]
  $$
  +\frac{\eX^{k'}_{(\hat l^{s*}_{r-1},\hat l^{s*}_r)} \ee1_{(\hat l^{s*}_{r-1},\hat l^{s*}_r)}}{(\hat l^{s*}_r-\hat l^{s*}_{r-1})^{1/2}}
  $$
  Using the claim (i), the assumptions (H2), (H3) we obtain that the last term converges to the sum of two normal distributions. Then, for $n \rightarrow \infty$: 
  $\eP[\exists k \in\{1, \cdots, p \}; k \not \in {\cal A}^*_{(l^0_{r-1},l^0_r)} \textrm{ and } k \in  {\cal A}^*_{n,(\hat l^{s*}_{r-1},\hat l^{s*}_r)} ]\rightarrow  0$.\\
  In the case $(\hat l^{s*}_{r-1},\hat l^{s*}_r) \not \subseteq (l^0_{r-1},l^0_r)$, suppose, without loss of generality, that $\hat l^{s*}_{r-1}  \leq l^0_{r-1} < \hat l^{s*}_r \leq  l^0_r$ (other cases are  similar). So, we have the decomposition: $\eX^{k'}_{(\hat l^{s*}_{r-1},\hat l^{s*}_r)} (\eY_{(\hat l^{s*}_{r-1},\hat l^{s*}_r)}- \eeX'_{(\hat l^{s*}_{r-1},\hat l^{s*}_r)} \hat \ef^{s*}_{(\hat l^{s*}_{r-1},\hat l^{s*}_r)})$\\$=\eX^{k'}_{(\hat l^{s*}_{r-1}, l^0_{r-1})} (\eY_{(\hat l^{s*}_{r-1},l^0_{r-1})}-\eeX'_{(\hat l^{s*}_{r-1},l^0_{r-1})} \hat \ef^{s*}_{(\hat l^{s*}_{r-1},\hat l^{s*}_r)})+\eX^{k'}_{(l^0_{r-1},\hat l^{s*}_r)}(\eY_{(l^0_{r-1},\hat l^{s*}_r)}-\eeX'_{(l^0_{r-1},\hat l^{s*}_r)} \hat \ef^{s*}_{(\hat l^{s*}_{r-1},\hat l^{s*}_r)})\equiv A_n+D_n$. As previously, $D_n/(\hat l^{s*}_r-\hat l^{s*}_{r-1} )^{1/2}$ converge to  the sum of two normal distributions. For $A_n$:
$A_n/(\hat l^{s*}_r-\hat l^{s*}_{r-1} )^{1/2}=\eX^{k'}_{(\hat l^{s*}_{r-1}, l^0_{r-1})}  \eeX_{(\hat l^{s*}_{r-1}, l^0_{r-1})}  (\ef^0_{r-1}-\hat \ef^{s*}(\hat l^{s*}_{r-1},\hat l^{s*}_r))  (\hat l^{s*}_r-\hat l^{s*}_{r-1} )^{1/2}+ \eX^{k'}_{(\hat l^{s*}_{r-1}, l^0_{r-1})} \cdot$ $\ee1_{(\hat l^{s*}_{r-1}, l^0_{r-1})}(\hat l^{s*}_r-\hat l^{s*}_{r-1} )^{1/2} \equiv A_{1n}+A_{2n}$. But $A_{2n}=o_{\eP}(1)$. Combining  (H2) and Theorem  \ref{Theorem 6} we have: $A_{1n}=\frac{\eX^{k'}_{(\hat l^{s*}_{r-1}, l^0_{r-1})} \cdot \eeX_{(\hat l^{s*}_{r-1}, l^0_{r-1})} \cdot (\ef^0_{r-1}-\hat \ef^{s*}(\hat l^{s*}_{r-1},\hat l^{s*}_r))}{\hat l^{s*}_r-\hat l^{s*}_{r-1}} \cdot \frac{\hat l^{s*}_r-\hat l^{s*}_{r-1}}{(\hat l^{s*}_r-\hat l^{s*}_{r-1})^{1/2}}=O_{\eP}(1)o_{\eP}(1)=o_{\eP}(1)$. Then $A_n=o_{\eP}(1)$. Consequently, there exists a constant $M >0$ such that $\lim_{n \rightarrow \infty} \eP[|A_n+D_n| <M]=0$. Hence $\eP[\exists k \in\{1, \cdots, p \}; k \not \in {\cal A}^*_{(l^0_{r-1},l^0_r)} \textrm{ and } k \in  {\cal A}^*_{n,(\hat l^{s*}_{r-1},\hat l^{s*}_r)} ]\rightarrow  0$,  for $n \rightarrow \infty$.
  
\hspace*{\fill}$\blacksquare$ \\
 
 \textbf{Acknowledgements:} The author would like to thank the Editor and the referee for their constructive comments and suggestions which helped to improve the quality of the paper.




\begin{thebibliography}{3}
\bibitem[Babu(1989)]{Babu:89}
Babu, G. J.,(1989).
\newblock Strong representations for LAD estimators in linear models.
 \newblock{\it Probability Theory and Related Fields}, \textbf{83},   547-558.
\bibitem[Bai(1998)]{Bai:98}
Bai, J.,(1998). 
\newblock Estimation of multiple-regime regressions with least absolute deviation.
\newblock{\it Journal of Statistical Planning Inference}, \textbf{74},  103-134.
\bibitem[Bai and Perron(1998)]{Bai:Perron:98} Bai,  J.,  Perron P.,(1998),
\newblock Estimating and testing linear models with multiple structural changes,
\newblock {\em Econometrica} {\bf 66}(1), 47-78.
\bibitem[Bickel et al.(2009)]{Bickel:Ritov:Tsybakov:09}
Bickel, P. J., Ritov, Y. and Tsybakov, A. B.,(2009). 
\newblock  Simultaneous analysis of lasso and Dantzig selector.
\newblock  {\it The Annals of Statistics}, \textbf{37}(4),  1705-1732.
\bibitem[Ciuperca(2009)]{Ciuperca:09}
Ciuperca G.,(2009).
\newblock The  M-estimation in a multi-phase random nonlinear model. 
\newblock {\it Statistics and  Probability Letters}, {\bf 75}(5), 573-580.
\bibitem[Ciuperca(2011a)]{Ciuperca:11a}
Ciuperca, G.,(2011a). 
\newblock Estimating nonlinear regression with and without change-points by the LAD-method.
\newblock{\it Annals of the Institute of Statistical Mathematics},   \textbf{63}(4), 717-743.
\bibitem[Ciuperca(2011b)]{Ciuperca:11b}
 Ciuperca, G.,(2011b). 
\newblock Penalized least absolute deviations estimation for  nonlinear model with change-points.
\newblock {\it Statistical Papers},  \textbf{52}(2), 371-390.
\bibitem[Fan and Li(2001)]{Fan:Li:01}
Fan, J. and  Li, R.,(2001). 
\newblock  Variable selection via nonconcave penalized likelihood and its oracle properties.
\newblock  {\it Journal of the American Statistical Association}, \textbf{96}(456),  1348-1360.
\bibitem[Foster et al.(2009)]{Foster:Verbyla:Pitchford:09}
Foster, S. D., Verbyla, A. P. and Pitchford, W. S.,(2009). 
\newblock Estimation, prediction and inference for the LASSO random effects model.
\newblock {\it The Australian $\&$ New Zealand Journal of Statistics }, \textbf{51}(1),  43-61.
 \bibitem[Harchaoui and L\'evy-Leduc(2010)]{Harchaoui:Levy:10}
 Harchaoui, Z. and L\'evy-Leduc, C.,(2010).
 \newblock Multiple change-point estimation with a total variation penalty.
 \newblock {\it Journal of the American Statistical Association}, {\bf 105}(492), 1480-1493.
\bibitem[Kim and Kim(2008)]{Kim:Kim:08}
Kim, J. and Kim, H. J.,(2008). 
\newblock Asymptotic results in segmented multiple regression.
\newblock{\it Journal of Multivariate Analysis}, \textbf{99}(9),  2016-2038.
\bibitem[Knight and Fu(2000)]{Knight:Fu:00}
Knight, K. and Fu, W.,(2000). 
\newblock Asymptotics for LASSO-type estimators.
\newblock {\it The Annals of Statistics}, \textbf{28}(5),  1356-1378.
\bibitem[P\"{o}tscher and Schneider(2009)]{Potscher:Schneider:09}
P\"{o}tscher, B. M. and Schneider, U.,(2009). 
\newblock On the distribution of the adaptive LASSO estimator.
\newblock {\it Journal of Statistical Planning Inference}, \textbf{139},  2775-2790.
 \bibitem[Koul and Qian(2002)]{Koul:Qian:02}Koul, H.L.,   Qian, L.,(2002), 
\newblock  Asymptotics of maximum likelihood estimator in a two-phase linear regression model.
\newblock {\em Journal of Statistical Planning and Inference}  {\bf 108}, 99-119.
\bibitem[Tibshirani(1996)]{Tibshirani:96}
Tibshirani, R.,(1996).
\newblock Regression shrinkage and selection via the LASSO.
\newblock{\it Journal of the Royal Statistical Society}, Ser. B, {\bf 58}, 267-288.
\bibitem[Wei et al.(2011)]{Wei:Huang:Li:11}
Wei, F., Huang, J. and Li, H.,(2011).
\newblock Variable selection and estimation in high-dimensional varying-coefficient models. 
\newblock {\it Statistica Sinica} to appear, \textbf{21}(4), doi:10.5705/ss.2009.316.  
\bibitem[Wu(2008)]{Wu:08}
Wu, Y.,(2008).
\newblock Simultaneous change point analysis and variable selection in a regression problem.
\newblock {\it Journal of Multivariate Analysis}, \textbf{99}(9),  2154-2171.
\bibitem[Xu and Ying(2010)]{Xu:Ying:10}
Xu, J. and Ying, Z.,(2010). 
\newblock Simultaneous estimation and variable selection in median regression using Lasso-type penalty.
\newblock {\it Annals of the Institute of Statistical Mathematics},  \textbf{62}, 487-514.
\bibitem[Yao(1988)]{Yao:88}
Yao, Y.C.,(1988), 
\newblock  Estimating the number of change-points via Schwarz's criterion.
\newblock {\em Statistics and Probability Letters}, {\bf 6},  181-189.
\bibitem[Zou(2006)]{Zou:06}
Zou, H.,(2006). 
\newblock  The adaptive Lasso and its oracle properties.
\newblock  {\it Journal of the American Statistical Association}, \textbf{101}(476),  1418-1428.
\end{thebibliography}
\end{document}